\documentclass[pdftex,11pt,reqno]{amsart}
\usepackage[utf8]{inputenc}
\usepackage{enumerate}
\usepackage{url}
\usepackage{cite}
\usepackage{comment}
\usepackage{xcolor}
\usepackage{amsrefs}
\usepackage{amsfonts, amsmath, amssymb, amscd, amsthm, bm, cancel}
\usepackage[
linktocpage=true,colorlinks,citecolor=magenta,linkcolor=blue,urlcolor=magenta]{hyperref}
\usepackage[margin=1in]{geometry}
 \newtheorem{thm}{Theorem}[section]
  
 \newtheorem{cor}[thm]{Corollary}
  
 \newtheorem{lem}[thm]{Lemma}
 
 \newtheorem{prop}[thm]{Proposition}
 \theoremstyle{definition}
 \newtheorem{defn}[thm]{Definition}
  \newtheorem*{ack}{Acknowledgments}
 \theoremstyle{remark}
 \newtheorem{rem}[thm]{Remark}

 \numberwithin{equation}{section}
 \allowdisplaybreaks

\renewcommand{\(}{\left(}
\renewcommand{\)}{\right)}
\renewcommand{\~}{\tilde}
\renewcommand{\-}{\overline}

\newcommand{\R}{\mathbb{R}}

\renewcommand{\H}{\mathbb{H}}

\renewcommand{\a}{\alpha}

\newcommand{\g}{\varphi}
\renewcommand{\d}{\delta}

\renewcommand{\k}{\kappa}
\renewcommand{\l}{\lambda}

\renewcommand{\t}{\theta}
\newcommand{\s}{\sigma}

\newcommand{\G}{\Gamma}

\renewcommand{\L}{\Lambda}

\newcommand{\ra}{\rightarrow}

\begin{document}
\title[Locally constrained curvature flows and geometric inequalities]{Locally constrained curvature flows and geometric inequalities in hyperbolic space}

\author{Yingxiang Hu}
\address{School of Mathematical Sciences, Beihang University, Beijing 100191, P.R. China}
\email{\href{mailto:huyingxiang@buaa.edu.cn}{huyingxiang@buaa.edu.cn}}

\author{Haizhong Li}
\address{Department of Mathematical Sciences, Tsinghua University, Beijing 100084, P.R. China}
\email{\href{mailto:lihz@tsinghua.edu.cn}{lihz@tsinghua.edu.cn}}

\author{Yong Wei}
\address{School of Mathematical Sciences, University of Science and Technology of China, Hefei 230026, P.R. China}
\email{\href{mailto:yongwei@ustc.edu.cn}{yongwei@ustc.edu.cn}}



\begin{abstract}
In this paper, we first study the locally constrained curvature flow of hypersurfaces in hyperbolic space, which was introduced by Brendle, Guan and Li \cite{BrendleGL}. This flow preserves the $m$th quermassintegral and decreases $(m+1)$th quermassintegral, so the convergence of the flow yields sharp Alexandrov-Fenchel type inequalities in hyperbolic space. Some special cases have been studied in \cite{BrendleGL}. In the first part of this paper, we show that h-convexity of the hypersurface is preserved along the flow and then the smooth convergence of the flow for h-convex hypersurfaces follows. We then apply this result to establish some new sharp geometric inequalities comparing the integral of $k$th Gauss-Bonnet curvature of a smooth h-convex hypersurface to its $m$th quermassintegral (for $0\leq m\leq 2k+1\leq n$), and comparing the weighted integral of $k$th mean curvature to its $m$th quermassintegral (for $0\leq m\leq k\leq n$). In particular, we give an affirmative answer to a conjecture proposed by Ge, Wang and Wu in 2015.

In the second part of this paper, we introduce a new locally constrained curvature flow using the shifted principal curvatures. This is natural in the context of h-convexity. We prove the smooth convergence to a geodesic sphere of the flow for h-convex hypersurfaces, and provide a new proof of the geometric inequalities proved by Andrews, Chen and the third author of this paper in 2018. We also prove a family of new sharp inequalities involving the weighted integral of $k$th shifted mean curvature for h-convex hypersurfaces, which as application implies a higher order analogue of Brendle, Hung and Wang's \cite{BHW2016} inequality.
\end{abstract}

\maketitle
\tableofcontents

\section{Introduction}\label{sec:1}
For any smooth convex body $\Omega$ in hyperbolic space $\mathbb{H}^{n+1}$ with boundary $M=\partial\Omega$ , the $k$th quermassintegral $W_k$ is defined  as the measure of the set of totally geodesic $k$-dimensional subspaces which intersect $\Omega$. $W_k$ can be computed as linear combination of integral of $k$th mean curvatures and of the enclosed volume (see \cite{Sant2004})\footnote{Note that this definition is different with the definition given in \cite{WX14} by a constant multiple $\frac{n+1-k}{n+1}$, see \cite{ACW2018}.}
\begin{align*}
& W_0(\Omega)= \mathrm{Vol}(\Omega),\qquad W_{1}(\Omega)= \frac{1}{n}|M|,\\
& W_{k+1}(\Omega)= \frac 1{n-k} \int_ME_k(\kappa)d\mu-\frac{k}{n-k}W_{k-1}(\Omega),\quad k=1,\cdots,n-1,
\end{align*}
where $E_k(\kappa)$ is defined as the normalized $k$th elementary symmetric function of the principal curvatures $\kappa=(\kappa_1,\cdots,\kappa_n)$ of $M$. The quermassintegrals satisfy the following nice variational property:
\begin{equation*}
  \frac d{dt}W_k(\Omega_t)=\int_{\partial\Omega_t}\eta E_k(\kappa)d\mu_t
\end{equation*}
along any normal variation with speed $\eta$.

\begin{defn}
A smooth bounded domain $\Omega$ in hyperbolic space $\mathbb{H}^{n+1}$ is said to be h-convex (resp. strictly h-convex) if the principal curvatures of the boundary $\partial\Omega$ satisfy $\kappa_i\geq 1$ (resp. $\kappa_i>1$) for all $i=1,\cdots,n$. Equivalently, $\Omega$ is h-convex if for any boundary point $p\in\partial\Omega$ there exists a horosphere enclosing $\Omega$ and touching $\Omega$ at $p$.
\end{defn}
In 2014, Wang and Xia \cite{WX14} proved the following Alexandrov-Fenchel type inequalities for h-convex domains.

\medskip\noindent
{\bf Theorem A} (\cite{WX14}).  
 \emph{Let $\Omega$ be a bounded smooth h-convex domain in $\mathbb H^{n+1}$. Then there holds
	\begin{align}\label{WX-inequality}
	W_k(\Omega) \geq f_k \circ f^{-1}_\ell(W_\ell(\Omega)),  \quad 0\leq \ell<k\leq n.
	\end{align}
	Equality holds if and only if $\Omega$ is a geodesic ball. Here $f_k:[0,\infty)\ra \R^{+}$ is a monotone function defined by $f_k(r)=W_k(B_r)$, the $k$th quermassintegral for the geodesic ball of radius $r$, and $f_\ell^{-1}$ is the inverse function of $f_\ell$. }

\medskip

\noindent The method to prove \eqref{WX-inequality} in \cite{WX14} is by using the quermassintegral preserving curvature flow $X:M^n\times [0,T)\to \mathbb{H}^{n+1}$:
\begin{equation*}
  \frac{\partial}{\partial t}X=\left(\phi(t)-\left(\frac{E_k}{E_\ell}\right)^{\frac 1{k-\ell}}\right)\nu,
\end{equation*}
where $\phi(t)$ is a global term that is chosen to preserve the $\ell$th quermassintegral of the enclosed domains of $M_t=X(M^n,t)$, and $\nu$ is the unit outward normal of $M_t$. The inequality \eqref{WX-inequality} with $k=3$, $\ell=1$ was proved earlier by the second and third authors of this paper with Xiong in \cite{LWX14} for star-shaped domains with $2$-convex boundary. Recently, the third author of this paper with Andrews and Chen \cite{ACW2018} proved \eqref{WX-inequality} with $k=1,\cdots,n$ and $\ell=0$ for domains with boundary having positive sectional curvatures; and the first two authors of this paper with Andrews \cite{AHL19} proved \eqref{WX-inequality} with $k=n-1$ and $\ell=n-1-2m ~(0<2m<n-1)$ for strictly convex domains.

In this paper, we view hyperbolic space as the warped product space $\H^{n+1}=[0,\infty)\times \mathbb S^{n}$ equipped with metric
\begin{equation*}
  \bar{g}=dr^2+\l(r)^2 g_{\mathbb S^n},
\end{equation*}
where $\l(r)=\sinh r$. A hypersurface $M$ in the hyperbolic space $\mathbb{H}^{n+1}$ is called star-shaped if its support function $u=\langle \lambda(r)\partial_r,\nu\rangle$ is positive everywhere on $M$. Equivalently, $M$ can be expressed as a graph of a smooth function over the unit sphere $\mathbb{S}^n$:
\begin{equation*}
  M=\{(r(\theta),\theta),~\theta\in \mathbb{S}^n,~r\in C^{\infty}(\mathbb{S}^n)\}.
\end{equation*}
Let $M_0$ be a smooth, star-shaped and $m$-convex hypersurface in $\mathbb{H}^{n+1}$, which is given by a smooth immersion $X_0:M^n \ra \mathbb{H}^{n+1}$. Brendle, Guan and Li \cite{BrendleGL} introduced the following locally constrained inverse curvature type flow $X:M \times [0,T) \ra \mathbb{H}^{n+1}$
\begin{equation}\label{s1:flow-BGL}
\left\{\begin{aligned}
\frac{\partial}{\partial t}X(x,t)=& \( \l'(r)\frac{E_{m-1}(\kappa)}{E_m(\kappa)}-u\) \nu(x,t),\quad m=1,\cdots,n,\\
X(\cdot,0)=& X_0(\cdot),
\end{aligned}\right.
\end{equation}
where $\lambda'(r)=\cosh r$, and $\kappa=(\kappa_1,\cdots,\kappa_n)$ are principal curvatures of $M_t=X(M,t)$.  The equation \eqref{s1:flow-BGL} is a pointwise defined parabolic PDE and preserves the $m$th quermassintegral of the enclosed domain $\Omega_t$ (by the Minkowski formula \eqref{s2:MF}), and thus is kind of locally constrained curvature flow. As the flow \eqref{s1:flow-BGL} decreases the $(m+1)$th quermassintegral $W_{m+1}(\Omega_t)$, this provides a potential method to prove the inequalities \eqref{WX-inequality} for starshaped and $m$-convex domains in hyperbolic space. The following result has been proved in \cite{BrendleGL} (see also a recent survey article \cite{GL19} by Guan and Li).

\medskip\noindent
{\bf Theorem B} (\cite{BrendleGL}).  
 \emph{Let $M_0$ be a smooth closed hypersurface in $\mathbb H^{n+1}$ satisfying either
 \begin{enumerate}
   \item $M_0$ is strictly convex, $m=n$; or
   \item $M_0$ is star-shaped and m-convex ($m=1,\cdots,n-1$), and a gradient bound
   \begin{equation}\label{s1:BGL-cond}
     \max_{x\in \mathbb{S}^n}|D\ln\cosh r|^2\leq 12+3\min_{x\in \mathbb{S}^n}\sinh^2 r
   \end{equation}
   is satisfied on $M_0$.
 \end{enumerate}
 Then the solution to the flow \eqref{s1:flow-BGL} exists for all time $t\in [0,\infty)$, and the solution $M_t$ converges exponentially to a geodesic sphere as $t\to\infty$ in $C^\infty$-topology. }

\medskip

\noindent The uniform $C^1$ estimate is a crucial step to prove the convergence of the flow \eqref{s1:flow-BGL}. For case (1), Brendle, Guan and Li proved that if $M_0$ is strictly convex, then the quotient $E_n(\kappa)/{E_{n-1}(\kappa)}$ is bounded from below by a positive constant and the strict convexity is preserved. Then $C^1$ estimate follows from the $C^0$ estimate. In general case, they proved that the condition \eqref{s1:BGL-cond} on the initial hypersurface is preserved along the flow \eqref{s1:flow-BGL} and this provides a uniform $C^1$ estimate of the solution. As applications of Theorem B, the inequality \eqref{WX-inequality} with $k=n$ and $\ell<n$ holds for strictly convex domains, and for general $0\leq \ell<k\leq n-1$, the inequality \eqref{WX-inequality} holds for star-shaped domains with boundary satisfying \eqref{s1:BGL-cond}.  Combining the convergence result in Theorem B with Gerhardt's \cite{Ge11} result on inverse mean curvature flow, Brendle, Guan and Li \cite{BrendleGL} also proved the inequality \eqref{WX-inequality} with $k=2$ and $\ell=1$ for star-shaped domains with mean convex boundary. 

We remark that the idea of using Minkowski (type) formula to define locally constrained curvature flows has been explored by Guan, Li and Wang \cite{GL15,GLW16} in space forms and certain warped product spaces. See \cite{LS19,Sche19,SWX18,SX19,WX19,WeiX19} for recent new results on locally constrained curvature flows.

\subsection{Brendle-Guan-Li's flow and applications}$\ $

In this paper, we first prove that if the initial hypersurface is h-convex, then the solution of the flow \eqref{s1:flow-BGL} remains to be h-convex for positive time. The h-convexity together with $C^0$ estimate of the solution immediately implies the uniform $C^1$ estimate. The curvature estimate and higher regularity estimate follow from the same steps in \cite{BrendleGL}. Then we obtain the following convergence result of the flow \eqref{s1:flow-BGL} for h-convex initial hypersurface in hyperbolic space.
\begin{thm}\label{main-thm-I}
Let $X_0: M^n \ra \mathbb \H^{n+1} (n\geq 2)$ be a smooth, closed and h-convex hypersurface in $\mathbb H^{n+1}$. Then the flow \eqref{s1:flow-BGL} has a unique smooth solution for all time $t\in [0,\infty)$, $M_t=X_t(M)$ is strictly h-convex for each $t>0$ and it converges smoothly and exponentially to a geodesic sphere of radius $r_\infty$ determined by $W_m(B_{r_\infty})=W_m(\Omega_0)$ as $t\to \infty$.
\end{thm}

As introduced in \cite{BrendleGL}, the smooth convergence of the flow \eqref{s1:flow-BGL} implies the Alexandrov-Fenchel inequalities \eqref{WX-inequality}. Therefore, Theorem \ref{main-thm-I} provides a new proof of the inequalities \eqref{WX-inequality} for h-convex domains in hyperbolic space.

Theorem \ref{main-thm-I} is also powerful to prove new inequalities. Let $(M,g)$ be a hypersurface in $\H^{n+1}$. The $k$th Gauss-Bonnet curvature $L_k$ of the metric $g$ is defined by
\begin{align}\label{Gauss-Bonnet-curvature}
L_k(g):=\frac{1}{2^k}\d_{j_1 j_2 \cdots j_{2k-1} j_{2k}}^{i_1 i_2 \cdots i_{2k-1} i_{2k}} R_{i_1 i_2}{}^{j_1 j_2} \cdots R_{i_{2k-1}i_{2k}}{}^{j_{2k-1}j_{2k}},
\end{align}
where $R_{ij}{}^{kl}$ is the Riemannian curvature tensor in the local coordinates with respect to the metric $g$, and the generalized Kronecker delta is defined by
\begin{align*}
\d_{j_1j_2\cdots j_r}^{i_1i_2\cdots i_r}=\det\(\begin{matrix}
\d^{i_1}_{j_1} & \d^{i_1}_{j_2} & \cdots & \d^{i_1}_{j_r} \\
\d^{i_2}_{j_1} & \d^{i_2}_{j_2} & \cdots & \d^{i_2}_{j_r} \\
\vdots        &    \vdots      & \vdots & \vdots \\
\d^{i_r}_{j_1} & \d^{i_r}_{j_2} & \cdots & \d^{i_r}_{j_r}
\end{matrix}\).
\end{align*}
The $k$th Gauss-Bonnet curvature $L_k$ of the induced metric of a hypersurface in hyperbolic space can be expressed in terms of the $k$th mean curvatures (see \cite[Lemma 3.1]{GeWW14})
\begin{align}\label{s1.Lk}
L_k(g)=\binom{n}{2k}(2k)!\sum_{j=0}^{k}(-1)^{j}\binom kj E_{2k-2j}(\kappa).
\end{align}
Applying Theorem \ref{main-thm-I}, we prove the following new geometric inequalities for h-convex hypersurfaces in hyperbolic space.
\begin{thm}\label{thm-geom-inequality}
	Let $M=\partial \Omega$ be a smooth, closed and h-convex hypersurface in $\mathbb H^{n+1}$. Then for any $0\leq m \leq 2k+1 \leq n$, there holds
	\begin{align}\label{GWW-inequality}
		\int_{M}L_k d\mu \geq g_k \circ f^{-1}_m(W_m(\Omega)),
	\end{align}
	where $L_k$ is the $k$th Gauss-Bonnet curvature of the induced metric on $M$ defined by \eqref{Gauss-Bonnet-curvature}, $g_k(r)=\binom{n}{2k}(2k)!\omega_{n}\sinh^{n-2k}r$ and $\omega_n$ is the area of the unit sphere $\mathbb S^n$ in $\mathbb R^{n+1}$. Equality holds in \eqref{GWW-inequality} if and only if $M$ is a geodesic sphere.
\end{thm}
The inequality \eqref{GWW-inequality} with $m=1$ was proved by Ge, Wang and Wu in \cite[Theorem 1.1]{GeWW14} for h-convex hypersurfaces, and later by the first two authors of this paper for hypersurfaces with nonnegative sectional curvatures \cite[Theorem 1]{Hu-Li2019}. The inequality \eqref{GWW-inequality} with $m=0$ follows from the case $m=1$ and isoperimetric inequality.

We also prove the following weighted Alexandrov-Fenchel inequalities for h-convex hypersurfaces in hyperbolic space.
\begin{thm}\label{thm-weighted-AF-inequality}
	Let $M=\partial \Omega$ be a smooth, closed and h-convex hypersurface in $\mathbb H^{n+1}$. For any $1\leq k \leq n$ and $0 \leq m\leq k$, there holds
\begin{align}\label{weighted-AF-inequality}
\int_{M} \lambda'(r) E_k(\kappa) d\mu\geq h_k \circ f_m^{-1}(W_m(\Omega)),
\end{align}
where $h_k(r)=\omega_{n}\cosh^{k+1}r\sinh^{n-k}r$. Equality holds in \eqref{weighted-AF-inequality} if and only if $M$ is a geodesic sphere centered at the origin.
\end{thm}

In particular, if $m=1$, we give an affirmative answer to a conjecture proposed by Ge, Wang and Wu \cite[Conjecture 9.1]{Ge-Wang-Wu2015}.
\begin{cor}\label{cor-weighted-AF-inequality}
	Let $M$ be a smooth, closed and h-convex hypersurface in $\mathbb H^{n+1}$. Then for any $1 \leq k \leq n$, there holds
	\begin{align}\label{weighted-AF-inequality-II}
	\int_{M} \lambda'(r) E_k(\kappa) d\mu \geq \omega_{n} \left(\(\frac{|M|}{\omega_n}\)^{\frac{2(n+1)}{n(k+1)}}+\(\frac{|M|}{\omega_n}\)^\frac{2(n-k)}{n(k+1)}\right)^\frac{k+1}{2}.
	\end{align}
	Equality holds in \eqref{weighted-AF-inequality-II} if and only if $M$ is a geodesic sphere centered at the origin.
\end{cor}
For $k=1$, \eqref{weighted-AF-inequality-II} was proved by de Lima and Girao \cite{deLima-Girao2016} for strictly mean-convex and star-shaped hypersurfaces. Later, for $k$ being odd, Ge, Wang and Wu \cite{Ge-Wang-Wu2015} proved \eqref{weighted-AF-inequality-II} for h-convex hypersurfaces. A partial answer to the conjecture \cite[Conjecture 9.1]{Ge-Wang-Wu2015} was also obtained recently by Girao, Pinheiro, Pinheiro, and Rodrigues \cite{Gir19}. For $k=n$, by the work of Brendle, Guan and Li \cite{BrendleGL}, Theorem \ref{thm-weighted-AF-inequality} and Corollary \ref{cor-weighted-AF-inequality} also hold for strictly convex hypersurfaces.

\subsection{New locally constrained curvature flow}$\ $

Recently, the third author of this paper with Andrews and Chen \cite{ACW2018} introduced a family of new quermassintegrals $\widetilde{W}_k(\Omega)$ of bounded domains in hyperbolic space $\mathbb{H}^{n+1}$
\begin{equation*}
  \widetilde{W}_k(\Omega)~:=~\sum_{i=0}^k(-1)^{k-i}\binom{k}{i}W_i(\Omega),\quad k=0,\cdots,n,
\end{equation*}
which are natural under the condition of h-convexity. These new quermassintegrals satisfy the following variational formula
\begin{equation}\label{s1:var2}
  \frac d{dt}\widetilde{W}_k(\Omega_t)=\int_{\partial\Omega_t}\eta E_k(\tilde{\kappa})d\mu_t
\end{equation}
along any variation in the direction of outward normal with speed function $\eta$, where $\tilde{\kappa}_i=\kappa_i-1$ are the shifted principal curvatures of the flow hypersurface $M_t=\partial\Omega_t$. It was proved that the new quermassintegrals satisfy the following Alexandrov-Fenchel type inequalities which improved Wang and Xia's inequalities \eqref{WX-inequality}.

\medskip\noindent
{\bf Theorem C} (\cite{ACW2018}).  
 \emph{Let $\Omega$ be a bounded smooth strictly h-convex domain in $\mathbb H^{n+1}$. Then there holds
\begin{equation}\label{s1:ACW}
   \widetilde{W}_{k}(\Omega)\geq \tilde{f}_{k}\circ \tilde{f}_\ell^{-1}(\widetilde{W}_\ell(\Omega)),\qquad 0 \leq \ell<k\leq n.
 \end{equation}
Equality holds if and only if $\Omega$ is a geodesic ball. Here $\tilde{f}_k:[0,\infty)\to \mathbb{R}^{+}$ is a monotone function defined by $\tilde{f}_k(r)=\widetilde{W}_k(B_r)$, and $\tilde{f}_\ell^{-1}$ is the inverse function of $\tilde{f}_\ell$.}
\medskip

In the second part of this paper, we introduce a new locally constrained curvature flow $X:M \times [0,T) \ra \mathbb{H}^{n+1}$
\begin{equation}\label{s1:flow1}
\frac{\partial}{\partial t}X(x,t)= \( (\l'(r)-u)\frac{E_{m-1}(\tilde{\kappa})}{E_m(\tilde{\kappa})}-u\) \nu(x,t),\quad m=1,\cdots,n,
\end{equation}
where $\tilde{\kappa}_i=\kappa_i-1$ are the shifted principal curvatures of $M_t=X(M,t)$, and $E_{m}(\tilde{\kappa})$ is the normalized $m$th elementary symmetric function of $\tilde{\kappa}$. For each $m=1,\cdots,n$, let
\begin{equation*}
  \Gamma_m^+=\{x\in \mathbb{R}^n:~E_i(x)>0,~i=1,\cdots,m\}
\end{equation*}
be the Garding cone. We prove the following smooth convergence of the flow \eqref{s1:flow1} for smooth h-convex hypersurfaces in hyperbolic space with shifted principal curvature $\tilde{\kappa}\in \Gamma_m^+$.
\begin{thm}\label{main-thm-II}
Let $X_0: M^n \ra \mathbb \H^{n+1} (n\geq 2)$ be a smooth embedding such that $M_0=X_0(M)$ is a smooth, h-convex hypersurface in $\mathbb H^{n+1}$ with $\tilde{\kappa}\in \Gamma_m^+$. Then the flow \eqref{s1:flow1} has a smooth solution for all time $t\in [0,\infty)$, and $M_t=X_t(M)$ is strictly h-convex for each $t>0$ and converges smoothly and exponentially to a geodesic sphere of radius $r_\infty$ such that $\widetilde{W}_m(B_{r_\infty})=\widetilde{W}_m(\Omega_0)$, where $\Omega_0$ is the domain enclosed by $M_0$.
\end{thm}

\begin{rem}
Since $\lambda'-u>0$, the condition $\tilde{\kappa}\in \Gamma_m^+$ on the initial hypersurface guarantees the short time existence of the flow \eqref{s1:flow1}, see e.g., \cite{HP99}. In particular, the conclusion of Theorem \ref{main-thm-II} holds for strictly h-convex initial hypersurfaces.
\end{rem}
By the variational formula \eqref{s1:var2} and shifted Minkowski formula \eqref{s2:shift-MF}, the flow \eqref{s1:flow1} preserves the $m$th new quermassintegral $\widetilde{W}_m$ of the enclosed domains, as we have
\begin{equation*}
 \frac d{dt}\widetilde{W}_m(\Omega_t)=\int_{\partial\Omega_t}\left((\lambda'(r)-u)E_{m-1}(\tilde{\kappa})-uE_{m}(\tilde{\kappa})\right)d\mu_t=0.
\end{equation*}
Moreover, we also have
\begin{align*}
 \frac d{dt}\widetilde{W}_{m+1}(\Omega_t)=&\int_{\partial\Omega_t}\left((\lambda'(r)-u)\frac{E_{m+1}(\tilde{\kappa})E_{m-1}(\tilde{\kappa})}{E_{m}(\tilde{\kappa})}-uE_{m+1}(\tilde{\kappa})\right)d\mu_t \\
  \leq  & \int_{\partial\Omega_t}\left((\lambda'(r)-u)E_{m}(\tilde{\kappa})-uE_{m+1}(\tilde{\kappa})\right)d\mu_t~\leq~0.
\end{align*}
Therefore, Theorem \ref{main-thm-II} implies that the inequalities \eqref{s1:ACW} hold for $h$-convex domains in hyperbolic space with boundary satisfying $\tilde{\kappa}\in \Gamma_k^+$.

We also use this new flow to prove the following new weighted Alexandrov-Fenchel inequalities for h-convex hypersurfaces.
\begin{thm}\label{thm-shift-weighted-AF-inequality}
	Let $M=\partial \Omega$ be a smooth, h-convex hypersurface in $\mathbb H^{n+1}$. Let $k=1,\cdots, n$. If $\tilde{\kappa}\in \Gamma_k^+$, then there holds
	\begin{align}\label{shift-weighted-AF-inequality}
\int_{M}  (\l'(r)-u) E_{k}(\~\k) d\mu\geq \~h_k \circ \~f_k^{-1}(\widetilde{W}_k(\Omega)),
	\end{align}
	where $\~h_k(r)=\omega_{n} (\cosh r-\sinh r)^{k+1}\sinh^{n-k}r$. If $k<n$ and $\tilde{\kappa}\in \Gamma_{k+1}^+$, then there holds
	\begin{align}\label{shift-weighted-AF-inequality-I}
\int_{M} (\l'(r)-u) E_{k}(\~\k) d\mu\geq	\~h_k \circ \~f_{k+1}^{-1}(\widetilde{W}_{k+1}(\Omega)).
	\end{align}
	Equality holds in \eqref{shift-weighted-AF-inequality} or \eqref{shift-weighted-AF-inequality-I} if and only if $M$ is a geodesic sphere centered at the origin.
\end{thm}

	In general, $\~h_k$ may not be a strictly increasing function of $r$. If $k\leq \frac{n-1}{2}$, then $\~h_k$ is strictly increasing in $r$ and by \eqref{s1:ACW} we obtain
\begin{cor}\label{corollary-shift-weighted-AF-inequality}
		Let $M=\partial \Omega$ be a smooth, h-convex hypersurface in $\mathbb H^{n+1}$. Let $1 \leq k \leq \frac{n-1}{2}$. If $\tilde{\kappa}\in \Gamma_k^+$, then there holds
	\begin{align}\label{shift-weighted-AF-inequality-III}
	\int_{M} (\l'(r)-u)E_{k}(\~\k)d\mu \geq \~h_k \circ \~f_\ell^{-1}(\widetilde{W}_\ell(\Omega)), \quad 0\leq \ell\leq k.
	\end{align}
	If $\tilde{\kappa}\in \Gamma_{k+1}^{+}$, then there holds
	\begin{align}\label{shift-weighted-AF-inequality-IV}
	\int_{M} (\l'(r)-u)E_{k}(\~\k)d\mu \geq \~h_k \circ \~f_\ell^{-1}(\widetilde{W}_\ell(\Omega)), \quad 0\leq \ell\leq k+1.
	\end{align}
	Equality holds in \eqref{shift-weighted-AF-inequality-III} or \eqref{shift-weighted-AF-inequality-IV} if and only if $M$ is a geodesic sphere centered at the origin.
\end{cor}

The inequalities on weighted curvature integrals with weight involving $\lambda'(r)$ and the support function $u$ have been studied previously by Brendle, Hung and Wang \cite{BHW2016} in 2014.  They proved a Minkowski-type inequality for strictly mean convex and star-shaped hypersurface $M$ in hyperbolic space $\mathbb H^{n+1}$ of the following form (See Theorem 1.2 in \cite{BHW2016})
\begin{equation}\label{s1.BHW}
  \int_M\left(\lambda'E_1(\kappa)-u\right)d\mu\geq \omega_n \(\frac{|M|}{\omega_n}\)^\frac{n-1}{n}.
\end{equation}
Equality holds in \eqref{s1.BHW} if and only if $M$ is a geodesic sphere centered at the origin. By the Minkowski identity \eqref{s2:MF} and the expression \eqref{s1.Lk} for the Gauss-Bonnet curvature $L_1$, we can reformulate \eqref{s1.BHW} as
\begin{align}\label{s1-BHW-Lk-version}
\int_{M} u L_1 d\mu \geq \binom{n}{2}2!\omega_n \(\frac{|M|}{\omega_n}\)^\frac{n-1}{n}.
\end{align}
Our inequalities in Corollary \ref{corollary-shift-weighted-AF-inequality} can be used to generalize the inequality \eqref{s1-BHW-Lk-version} to higher order case.  For this purpose, we show that the $k$th Gauss-Bonnet curvature $L_k$ of the induced metric of a hypersurface in hyperbolic space can also be expressed in terms of a linear combination of the shifted $m$th mean curvatures with $k\leq m\leq 2k$ (see Lemma \ref{s9:lem-Lk} below)
\begin{align}\label{s1:Lk-shifted-curvature-integral}
L_k= \binom{n}{2k}(2k)! \sum_{i=0}^{k}\binom{k}{i}2^i E_{2k-i}(\~\kappa).
\end{align}
Applying Corollary \ref{corollary-shift-weighted-AF-inequality} and \eqref{s1:Lk-shifted-curvature-integral}, we have
\begin{cor}\label{corollary-weighted-Lk-inequality}
	Let $M=\partial \Omega$ be a smooth and h-convex hypersurface in $\mathbb H^{n+1}$. Let $2\leq k \leq \frac{n-1}{4}$. If $\tilde{\kappa}\in \Gamma_{2k-1}^+$, then there holds
	\begin{align}\label{shifed-weighted-Lk-inequality}
	\int_{M} u L_k d\mu \geq \~g_k \circ \~f_{\ell}^{-1}(\widetilde{W}_\ell(\Omega)),  \quad 0 \leq \ell < k,
	\end{align}
    where $\~g_k(r)=\binom{n}{2k}(2k)!\omega_n \sinh^{n+1-2k} r$. In particular, we have
    \begin{align}\label{shifed-weighted-BHW-inequality}
	\int_{M} u L_k d\mu \geq \binom{n}{2k}(2k)!\omega_n \(\frac{|M|}{\omega_n}\)^\frac{n+1-2k}{n}.
	\end{align}
	Equality holds in \eqref{shifed-weighted-Lk-inequality} or \eqref{shifed-weighted-BHW-inequality} if and only if $M$ is a geodesic sphere centered at the origin.
\end{cor}

The paper is organized as following. In Section \ref{sec:2}, we collect some basic properties on elementary symmetric functions, and recall some geometry of hypersurfaces in hyperbolic space. In Section \ref{sec:evl}, we derive the evolution equations along the flows \eqref{s1:flow-BGL} and \eqref{s1:flow1}.

In Section \ref{sec:h-conv1}, we apply the tensor maximum principle to show that the flow \eqref{s1:flow-BGL} preserves h-convexity of the evolving hypersurfaces, and then complete the proof of Theorem \ref{main-thm-I}. Then the convergence of the flow \eqref{s1:flow-BGL} will be used in Section \ref{sec:Ineqn} to prove the new geometric inequalities in Theorem \ref{thm-geom-inequality} and Theorem \ref{thm-weighted-AF-inequality}.

In Sections \ref{sec:new flow} -- \ref{sec:conver}, we study the new locally constrained curvature flow \eqref{s1:flow1} and prove the smooth convergence of the flow to a geodesic sphere. The new inequalities in Theorem \ref{thm-shift-weighted-AF-inequality},  Corollary \ref{corollary-shift-weighted-AF-inequality} and Corollary \ref{corollary-weighted-Lk-inequality} will be proved in Section \ref{sec:9}.

\begin{ack}
The authors would like to thank Professor Ben Andrews and Professor Pengfei Guan for helpful discussions, and the referees for their valuable comments and suggestions. The first author was supported by China Postdoctoral Science Foundation (No.2018M641317). The second author was supported by NSFC grant No.11831005, 11671224  and NSFC-FWO grant No.1196131001. The third author was supported by Discovery Early Career Researcher Award DE190100147 of the Australian Research Council and a research grant from University of Science and Technology of China.
\end{ack}

\section{Preliminaries}\label{sec:2}

In this section, we collect some preliminaries on elementary symmetric functions and geometry of hypersurfaces in hyperbolic space.
\subsection{Elementary symmetric functions}$\ $

We first review some properties of elementary symmetric functions. See \cite{Guan12} for more details.

For integer $m=1,\cdots,n$ and a point $\kappa=(\kappa_1,\cdots,\kappa_n)\in \mathbb{R}^n$, the (normalized) $m$th elementary symmetric function $E_m$ is defined by
\begin{align}\label{s2:Em-def}
E_m(\k)=\binom{n}{m}^{-1}\s_m(\k)=\binom{n}{m}^{-1}\sum_{1\leq i_1<\cdots<i_m \leq n}\k_{i_1}\cdots \k_{i_m}.
\end{align}
It is convenient to set $E_0(\kappa)=1$ and $E_m(\kappa)=0$ for $m>n$. The definition can be extended to symmetric matrices. Let $A\in \mathrm{Sym}(n)$ be an $n\times n$ symmetric matrix. Denote $\kappa=\kappa(A)$ the eigenvalues of $A$. Set $E_m(A)=E_m(\kappa(A))$. We have
\begin{align*}
  E_m(A) =& \frac{(n-m)!}{n!}\delta_{i_1,\cdots,i_k}^{j_1,\cdots,j_k}A_{i_1j_1}\cdots A_{i_kj_k},\quad m=1,\cdots,n.
\end{align*}
\begin{lem}\label{lem-Newton-tensor} Denote $\dot{E}_m^{ij}=\frac{\partial E_m}{\partial A_{ij}}$. Then we have
	\begin{align}
	\sum_{i,j}\dot{E}_m^{ij}A_{ij}= &m E_m,\label{newton-formula-1}\\
	\sum_{i,j}\dot{E}_m^{ij} \delta_i^j =& m E_{m-1},\label{newton-formula-2}\\
	\sum_{i,j}\dot{E}_m^{ij}(A^2)_{ij} =& n E_1 E_{m} -(n-m)E_{m+1},\label{newton-formula-3}
	\end{align}
where $(A^2)_{ij}=\sum_{k=1}^nA_{ik}A_{kj}$.
\end{lem}

\begin{lem}
If $\kappa\in \Gamma_m^+$, we have the following Newton-MacLaurin inequality
\begin{align}
  &E_{m+1}(\k) E_{k-1}(\k)\leq  E_k(\k) E_{m}(\k),\quad 1\leq k\leq m. \label{s2:Newt-2}
\end{align}
Equality holds if and only if $\kappa_1=\cdots=\kappa_n$.
\end{lem}
Using \eqref{newton-formula-2} -- \eqref{s2:Newt-2}, we have the following corollary.
\begin{cor}
Let $F(A)=E_m(A)/{E_{m-1}(A)}$ and $\k(A)\in \Gamma_{m}^{+}$. Then
\begin{align}
 & 1~\leq ~\sum_{i,j}\dot{F}^{ij}\delta_{i}^j \leq  ~ m \label{s2:0-1}\\
  & F^2~\leq~\sum_{i,j}\dot{F}^{ij}(A^2)_{ij} \leq ~ (n-m+1)F^2. \label{s2:0-2}
\end{align}
\end{cor}
%
\subsection{Hypersurfaces in hyperbolic space}$\ $

In this paper, the hyperbolic space $\H^{n+1}$ is viewed as a warped product manifold $\R^{+}\times \mathbb{S}^{n}$ equipped with the metric
\begin{align*}
\-g=dr^2+\l(r)^2 g_{\mathbb{S}^{n}},
\end{align*}
where $\l(r)=\sinh r$ and $g_{\mathbb{S}^{n}}$ is the standard metric on the unit sphere $\mathbb{S}^{n}\subset \mathbb R^{n+1}$. We define
\begin{equation*}
  \Phi(r)=\int_0^r \l(s)ds=\cosh r-1.
\end{equation*}
The vector field $V=\-\nabla \Phi=\l\partial_{r}$ on $\H^{n+1}$ is a conformal Killing field, i.e., $\-\nabla (\l \partial_r) =\l' \-g$.

Let $M$ be a closed hypersurface in $\H^{n+1}$ with unit outward normal $\nu$. The second fundamental form $h$ of $M$ is defined by $h(X,Y)=\langle \-\nabla_X \nu, Y \rangle$ for any tangent vectors $X,Y$ on $M$. The principal curvature $\k=(\k_1,\cdots,\k_n)$ are the eigenvalues of $h$ with respect to the induced metric $g$ on $M$. Let $\{x^1,\cdots,x^n\}$ be a local coordinate system of $M$ and denote $g_{ij}=g(\partial_i,\partial_j)$ and $h_{ij}=h(\partial_i,\partial_j)$. The Weingarten matrix is denoted by $\mathcal{W}=(h_i^j)$, where $h_i^j=g^{kj}h_{ik}$ and $(g^{ij})$ is the inverse matrix of $(g_{ij})$. Then the principal curvatures of $M$ are eigenvalues of the Weingarten matrix $\mathcal{W}$.

The following formulas hold for smooth hypersurfaces in $\H^{n+1}$, see e.g. \cite[Lems. 2.2 \& 2.6]{GL15}.

\begin{lem}\label{lem-1} Let $(M,g)$ be a smooth hypersurface in $\H^{n+1}$. Then $\Phi|_{M}$ satisfies
	\begin{align}\label{2.1}
	\nabla_i \Phi= \langle V,e_i\rangle,\quad \nabla_j\nabla_i \Phi=\l' g_{ij} -u h_{ij},
	\end{align}
	and the support function $u=\langle V,\nu \rangle$ satisfies
	\begin{align}\label{2.2}
	\nabla_i u=  \langle V,e_k\rangle h_i^k, \quad \nabla_j\nabla_i u= \langle V,\nabla h_{ij}\rangle+\l' h_{ij}-u (h^2)_{ij},
	\end{align}
	where $\{e_1,\cdots,e_n\}$ is a basis of the tangent space of $M$.
\end{lem}
We also have the well-known Minkowski formulas, see e.g. \cite[Proposition 2.5]{GL15}.
\begin{lem}
Let $M$ be a smooth closed hypersurface in hyperbolic space $\H^{n+1}$. Then
\begin{align}\label{s2:MF}
\int_{M} \l' E_m(\kappa)d\mu = \int_{M} uE_{m+1}(\kappa)d\mu,  \qquad  m=0,1,\cdots,n-1,
\end{align}
where $E_m(\kappa)$ is the $m$th mean curvature of $M$.
\end{lem}
From \eqref{s2:MF} we can derive the following Minkowski type formula for closed hypersurfaces in hyperbolic space.
\begin{lem}
Let $M$ be a closed hypersurface in $\mathbb{H}^{n+1}$. Let $\tilde{\kappa}_i=\kappa_i-1$ be the shifted principal curvatures of $M$. Then
\begin{equation}\label{s2:shift-MF}
  \int_M(\lambda'-u)E_m(\tilde{\kappa})d\mu=\int_M u E_{m+1}(\tilde{\kappa})d\mu,\qquad m=0,1,\cdots,n-1.
\end{equation}
\end{lem}
\proof
We have
\begin{align*}
  E_m(\tilde{\kappa}) =& \sum_{i=0}^m(-1)^{m-i}\binom mi E_i(\kappa) \\
  = & (-1)^m+\sum_{i=0}^{m-1}(-1)^{m-1-i}\binom m{i+1}E_{i+1}(\kappa),
\end{align*}
and
\begin{align*}
  E_{m+1}(\tilde{\kappa}) =& \sum_{i=0}^{m+1}(-1)^{m+1-i}\binom {m+1}i E_i(\kappa) \\
  = & (-1)^{m+1}+\sum_{i=0}^{m}(-1)^{m-i}\binom {m+1}{i+1}E_{i+1}(\kappa).
\end{align*}
Then
\begin{align*}
   & (\lambda'-u)E_m(\tilde{\kappa})-u E_{m+1}(\tilde{\kappa}) \\
  = & \lambda' E_m(\tilde{\kappa})-u(E_m(\tilde{\kappa})+E_{m+1}(\tilde{\kappa}))\\
  =&\lambda' E_m(\tilde{\kappa})-u(E_{m+1}(\kappa)+\sum_{i=0}^{m-1}(-1)^{m-i}(\binom {m+1}{i+1}-\binom m{i+1})E_{i+1}(\kappa))\\
  =& \lambda' E_m(\tilde{\kappa})-uE_{m+1}(\kappa)-u \sum_{i=0}^{m-1}(-1)^{m-i}\binom mi E_{i+1}(\kappa)\\
  =& \sum_{i=0}^m(-1)^{m-i}\binom mi (\lambda'E_i(\kappa)-uE_{i+1}(\kappa)).
\end{align*}
Integrating the above equation on the hypersurface and using \eqref{s2:MF} yields the new formula \eqref{s2:shift-MF}.
\endproof

\subsection{Parametrization by radial graph}\label{sec:2-3}$\ $

A smooth closed hypersurface $M$ in hyperbolic space is called {\em star-shaped} if the support function $u=\langle \l \partial_r,\nu \rangle>0$ everywhere on $M$. This is equivalent to that $M$ can be expressed as a radial graph in spherical coordinates $(r(\theta),\theta)$ in $\mathbb{H}^{n+1}$:
\begin{align*}
M= \left\{ (r(\theta),\theta) \in \R^{+}\times \mathbb S^n ~|~ \t \in \mathbb{S}^{n} \right\}.
\end{align*}
Let $\theta=(\theta^1,\cdots,\theta^n)$ be a local coordinate of $\mathbb{S}^n$. We write $\partial_i=\partial_{\theta^i}$ and $r_i=D_{i}r$, where $D$ denotes the Levi-Civita connection on $\mathbb{S}^n$ with respect to $g_{\mathbb{S}^n}$. Then the tangent space of $M$ is spanned by $\{e_i=\partial_i + r_i \partial_r, i=1,\cdots,n\}$. The induced metric on $M$ can be expressed as
\begin{align*}
g_{ij}=\l^2(r) e_{ij}+r_i r_j,
\end{align*}
where $e_{ij}=g_{\mathbb{S}^n}(\partial_i,\partial_j)$. The inverse matrix $(g_{ij})^{-1}$ is given by
\begin{align*}
g^{ij}= \frac 1{\l(r)^{2}} \(e^{ij}-\frac{r^i r^j}{\l(r)^2+|D r|^2}\),
\end{align*}
where $r^i=e^{ij}r_j$. For the convenience of the notations, we introduce a new function $\varphi:\mathbb S^{n} \ra \mathbb R$ by
$$
\varphi(\t)=\Psi(r(\t)),
$$
where $\Psi(r)$ is a positive function satisfying $\Psi'(r)=1/\l(r)$. Then we have $\varphi_i={r_i}/{\l(r)}$. Set
\begin{equation}\label{s2:v-def}
 v=\sqrt{1+|D \varphi|^2}.
\end{equation}
Then
\begin{equation*}
  g_{ij}= \l^2 (e_{ij}+\varphi_i \varphi_j), \qquad g^{ij}=\frac 1{\l^{2}}\(e^{ij}-\frac{\varphi^i\varphi^j}{v^2}\).
\end{equation*}
The unit outward normal vector on $M$ is given by
\begin{align*}
\nu= \frac 1v\left(\partial_r-\frac{r^i\partial_i}{\lambda(r)^2} \right).
\end{align*}
It follows that the support function
\begin{align}\label{s2:v-1}
u=\langle \l \partial_r,\nu \rangle =\frac{\l}{v}.
\end{align}
The second fundamental form $h$ is given by
\begin{align*}
h_{ij}=\frac{1}{v} \( -\lambda\varphi_{ij}+\frac{\lambda'}{\lambda}g_{ij}\),
\end{align*}
where the covariant derivatives are taken with respect to $g_{\mathbb{S}^n}$ on $\mathbb{S}^n$. The Weingarten matrix can be expressed as follows
\begin{align}\label{s2:h}
h_i^j=&g^{jk}h_{ik}=-\frac{1}{\l v}(e^{jk}-\frac{\varphi^j\varphi^k}{v^2})\varphi_{ik}+\frac{\lambda'}{\lambda v}\delta_i^j
\end{align}

It's well known that if a family of star-shaped hypersurfaces $M_t=X(M,t)$ satisfy the flow equation
\begin{equation}\label{s2:3-1}
  \frac{\partial}{\partial t}X=\eta\nu,
\end{equation}
where $\eta=\eta(x,t)$ is a smooth function on $M_t$, then we can express \eqref{s2:3-1} as the initial value problem of radial graphs $r=r(X(\t,t),t)$ over $\mathbb{S}^n$ in $\H^{n+1}$, where $(\t_1,\cdots,\t_n)$ is a local coordinates of $\mathbb S^n$:
\begin{equation*}
\left\{ \begin{aligned}
\frac{\partial}{\partial t} r =& \eta v, \quad \text{for $(\t,t)\in \mathbb{S}^n\times [0,\infty)$},\\
r(\cdot, 0)=& r_0(\cdot),
\end{aligned}\right.
\end{equation*}
which is also equivalent to the initial value problem for $\varphi$ on $\mathbb{S}^n$:
\begin{equation}\label{s2:gamma-evl}
\left\{ \begin{aligned}
\frac{\partial}{\partial t} \varphi =& \eta \frac{v}{\l}, \quad \text{for $(\t,t)\in \mathbb{S}^n\times [0,\infty)$},\\
\varphi(\cdot, 0)=& \varphi_0(\cdot),
\end{aligned}\right.
\end{equation}
where $v$ is the function defined in \eqref{s2:v-def}.

\section{Evolution equations}\label{sec:evl}

In this section, we derive the evolution equations along the flows \eqref{s1:flow-BGL} and \eqref{s1:flow1}. Denote
\begin{equation*}
 F(\kappa)=\frac{E_m(\kappa)}{E_{m-1}(\kappa)}.
\end{equation*}
Then $F$ is a smooth symmetric function on $\mathbb{R}^n$. Since the principal curvatures $\kappa=(\kappa_1,\cdots,\kappa_n)$ are eigenvalues of the Weingarten matrix $\mathcal{W}$,  we can view $F(\kappa)$ as a smooth symmetric function $F(\mathcal{W})$ of diagonalizable Weingarten matrix. Equivalently, we can view $F=F(h_{ij},g_{ij})$ as a smooth $GL(n)$-invariance function of the second fundamental form $h_{ij}$ and the metric $g_{ij}$. The invariance of $F$ implies that  $F=F(h_{ij},g_{ij})=F(g^{-1/2}hg^{-1/2},I)$. Thus $F$ is reduced to an $O(n)$-invariant function of the first argument. We write $\dot{F}^{kl}$ and $\ddot{F}^{kl,pq}$ as  the first and second derivatives of $F$  with respect to the first argument, and
\begin{equation*}
  \dot{F}^i(\kappa)=\frac{\partial F(\kappa)}{\partial \kappa_i},\qquad \ddot{F}^{ij}(\kappa)=\frac{\partial^2F(\kappa)}{\partial \kappa_i\partial \kappa_j}
\end{equation*}
for the derivatives of $F$ with respect to the principal curvatures. In a local orthonormal frame, we have (see e.g., \cite[Theorem 5.1]{And07})
\begin{equation*}
  \dot{F}^{kl}=\dot{F}^k(\kappa)\delta_{kl}
\end{equation*}
and
\begin{equation*}
  \ddot{F}^{kl,pq}B_{kl}B_{pq}=\ddot{F}^{kl}(\kappa)B_{kk}B_{ll}+2\sum_{k<l}\frac{\dot{F}^k(\kappa)-\dot{F}^l(\kappa)}{\kappa_k-\kappa_l}(B_{kl})^2
\end{equation*}
for any symmetric matrix $B$. The latter identity makes sense as a limit if $\kappa_k=\kappa_l$.

Now we derive the evolution equations along the flow \eqref{s1:flow-BGL}.
\begin{lem}
	Along the flow \eqref{s1:flow-BGL}, the induced metric $g=(g_{ij})$ and second fundamental form $(h_{ij})$ satisfy the following equations	
	\begin{align}
\frac{\partial}{\partial t}g_{ij}=&2\left(\frac{\lambda'(r)}{F}-u\right)h_{ij}\label{s3:1-gij}\\
	\frac{\partial}{\partial t}h_{ij}=&\frac{\lambda'}{F^2}\dot{F}^{kl}\nabla_k\nabla_lh_{ij}+\frac{\lambda'}{F^2}\ddot{F}^{kl,pq}\nabla_ih_{kl}\nabla_jh_{pq} \nonumber\\
	         &+\frac 1F\(\nabla_jF\nabla_i(\frac{\lambda'}F)+\nabla_iF\nabla_j(\frac{\lambda'}F)\)+\langle V,\nabla h_{ij}\rangle \nonumber\\
&+\(\frac{u}{F}+ \l'+\frac{\l'}{F^2}\dot{F}^{kl}((h^2)_{kl}+g_{kl})\) h_{ij} -2u(h^2)_{ij}-\(\frac{\l'}{F}+u\)g_{ij},\label{s3:1-hij}
	\end{align}
where $(h^2)_{ij}=\sum_k h_i^k h_{kj}$ and $V=\lambda(r)\partial_r$.
\end{lem}
\begin{proof}
The first equation \eqref{s3:1-gij} is well-known, see e.g., \cite{HP99}. For the second equation \eqref{s3:1-hij}, recall that along a general flow
\begin{align*}
\frac{\partial}{\partial t} X(x,t)= \eta(x,t)\nu(x,t)
\end{align*}
of smooth hypersurfaces in hyperbolic space, we have the following evolution equation on the second fundamental form of $M_t$ (see \cite{And94}):
\begin{align}\label{2.4}
\frac{\partial}{\partial t} h_{ij}=-\nabla_j \nabla_i \eta + \eta ((h^2)_{ij}+g_{ij}).
\end{align}
Let
\begin{equation*}
  \eta(x,t)=\frac{\lambda'(r)}{F(\kappa)}-u
\end{equation*}
in \eqref{2.4}. We have
\begin{align}
	\frac{\partial}{\partial t} h_{ij} = &-\nabla_j\nabla_i\(\frac{\lambda'}F-u\) +\(\frac{\lambda'}F-u\)\((h^2)_{ij}+g_{ij}\) \nonumber\\
	= &\frac{\lambda'}{F^2}\nabla_j\nabla_iF - \frac 1F\nabla_j\nabla_i\l'+\frac 1F\(\nabla_jF\nabla_i(\frac{\lambda'}F)+\nabla_iF\nabla_j(\frac{\lambda'}F)\)\nonumber\\
&+\nabla_j\nabla_iu+\(\frac{\lambda'}F-u\)\((h^2)_{ij}+g_{ij}\) . \label{s3:1-h-evl}
	\end{align}
By the Simons' identity (see e.g., \cite[Equation (2-7)]{And94}),
\begin{align}\label{s3:Sims}
  \nabla_j\nabla_iF= & \nabla_j(\dot{F}^{kl}\nabla_ih_{kl}) \nonumber\\
   =&\dot{F}^{kl}\nabla_j\nabla_ih_{kl}+\ddot{F}^{kl,pq}\nabla_ih_{kl}\nabla_jh_{pq}\nonumber\\
   =& \dot{F}^{kl}\nabla_k\nabla_lh_{ij} + \dot{F}^{kl}((h^2)_{kl}+g_{kl})h_{ij}-\dot{F}^{kl}h_{kl}((h^2)_{ij}+g_{ij})\nonumber\\
   &\quad +\ddot{F}^{kl,pq}\nabla_ih_{kl}\nabla_jh_{pq}.
\end{align}
Since $F$ is homogeneous of degree $1$, we have $F=\dot{F}^{kl}h_{kl}$. Substituting \eqref{2.1}, \eqref{2.2} and \eqref{s3:Sims} into \eqref{s3:1-h-evl}, we have
\begin{align*}
	\frac{\partial}{\partial t} h_{ij}= & \frac{\lambda'}{F^2}\dot{F}^{kl}\nabla_k\nabla_lh_{ij}+\frac{\lambda'}{F^2}\ddot{F}^{kl,pq}\nabla_ih_{kl}\nabla_jh_{pq}+\frac 1F\(\nabla_jF\nabla_i(\frac{\lambda'}F)+\nabla_iF\nabla_j(\frac{\lambda'}F)\)\nonumber\\
&+\frac{\lambda'}{F^2}\dot{F}^{kl}((h^2)_{kl}+g_{kl})h_{ij}-\frac{\lambda'}{F}((h^2)_{ij}+g_{ij})+\frac{1}F(uh_{ij}-\lambda'g_{ij})\nonumber\\
 &+\langle V,\nabla h_{ij} \rangle +\l'  h_{ij} -u (h^2)_{ij} +\(\frac{\lambda'}F-u\)\((h^2)_{ij}+g_{ij}\)\nonumber\\
= & \frac{\lambda'}{F^2}\dot{F}^{kl}\nabla_k\nabla_lh_{ij}+\frac{\lambda'}{F^2}\ddot{F}^{kl,pq}\nabla_ih_{kl}\nabla_jh_{pq}\nonumber\\
  &+\frac 1F\(\nabla_jF\nabla_i(\frac{\lambda'}F)+\nabla_iF\nabla_j(\frac{\lambda'}F)\)+\langle V,\nabla h_{ij} \rangle\nonumber\\
&+\(\frac{u}{F}+ \l'+\frac{\l'}{F^2}\dot{F}^{kl}((h^2)_{kl}+g_{kl})\) h_{ij} -2u(h^2)_{ij}-\(\frac{\l'}{F}+u\)g_{ij}.
\end{align*}
\end{proof}

To derive the evolution equations along the flow \eqref{s1:flow1}, we denote
\begin{equation*}
  F(\tilde{\kappa})=\frac{E_m(\tilde{\kappa})}{E_{m-1}(\tilde{\kappa})}
\end{equation*}
for the shifted principal curvatures $\tilde{\kappa}=(\kappa_1-1,\cdots,\kappa_n-1)$. Then we can also view $F$ as a smooth symmetric function of the shifted Weingarten matrix $S_i^j=h_i^j-\delta_i^j$.
\begin{lem}\label{s3:lem2}
	Along the flow \eqref{s1:flow1}, we have the following evolution equations.
	\begin{enumerate}[(i)]
\item The induced metric evolves as
\begin{equation}\label{s3:2-gij}
  \frac{\partial}{\partial t}g_{ij}=2\left(\frac{\lambda'(r)-u}{F(\tilde{\kappa})}-u\right)h_{ij}.
\end{equation}
		\item The shifted second fundamental form $S_{ij}=h_{ij}-g_{ij}$ evolves as
	\begin{align}\label{s3:2-hij}
	\frac{\partial}{\partial t} S_{ij} =&\frac{\lambda'-u}{F^2}\dot{F}^{kl}\nabla_k\nabla_lS_{ij}+\frac{\lambda'-u}{F^2}\ddot{F}^{kl,pq}\nabla_iS_{kl}\nabla_jS_{pq}\nonumber\\
	&+\frac 1F \(\nabla_jF\nabla_i(\frac{\lambda'-u}F)+\nabla_iF\nabla_j(\frac{\lambda'-u}F)\)\nonumber\\
&+\frac {1+F}F\langle V,\nabla S_{ij}\rangle+\left(\frac{\lambda'-u}{F^2}(\dot{F}^{kl}(S^2)_{kl}+F)+\lambda'-2u\right)S_{ij}\nonumber\\
&-\left(\frac{\lambda'-u}{F^2}\dot{F}^{kl}g_{kl}+\frac{u}F+2u\right)(S^2)_{ij}+(\lambda' -u) \left(1+\frac{\dot{F}^{kl}(S^2)_{kl}}{F^2}\right)g_{ij},
	\end{align}
	where $(S^2)_{ij}=\sum_kS_{i}^{k}S_{kj}$.
	  \item The curvature function $F$ evolves as
	\begin{align}\label{evol-F}
	\frac{\partial}{\partial t}F=& \frac{\lambda'-u}{F^2}\dot{F}^{ij}\nabla_i\nabla_jF +\frac 2F\dot{F}^{ij}\nabla_jF\nabla_i(\frac{\lambda'-u}F)+\frac {1+F}F\langle V,\nabla F\rangle\nonumber\\
 &-\frac{\lambda'}F\dot{F}^{ij}(S^2)_{ij}+\l' F +(\lambda'-u) (\dot{F}^{ij}g_{ij}-1).
	\end{align}
	\item The support function $u$ satisfies
\begin{align}\label{s3:u-evl}
  \frac{\partial}{\partial t}u =&\frac{\lambda'-u}{F^2}\dot{F}^{ij}\nabla_i\nabla_ju+\frac{1+F}F\langle V,\nabla u\rangle -\frac 1F\langle V,\nabla\lambda'\rangle\nonumber\\
  &+\frac{u(\lambda'-u)}{F^2}\dot{F}^{ij}(S^2)_{ij}-\frac{(\lambda'-u)^2}{F^2}\dot{F}^{ij}g_{ij}-\lambda'u+2u\frac{\lambda'-u}{F}.
\end{align}
\item The function $\lambda'(r)$ satisfies
\begin{align}\label{s3:lambd'-evl}
   \frac{\partial}{\partial t}\lambda'(r)=&\frac{\lambda'-u}{F^2}\dot{F}^{ij}\nabla_i\nabla_j\lambda'(r)+2u\frac{\lambda'-u}F-u^2-\frac{(\lambda'-u)^2}{F^2}\dot{F}^{ij}g_{ij}.
\end{align}
\end{enumerate}
\end{lem}
\begin{proof}
	(i) The evolution equation for the metric $g_{ij}$ is well known.

(ii) Let
\begin{equation*}
  \eta(x,t)=\frac{\lambda'(r)-u}{F(\tilde{\kappa})}-u
\end{equation*}
in the general formula \eqref{2.4}. We have
	\begin{align}
	\frac{\partial}{\partial t} h_{ij} = &-\nabla_j\nabla_i\(\frac{\lambda'-u}F-u\) +\(\frac{\lambda'-u}F-u\)\((h^2)_{ij}+g_{ij}\) \nonumber\\
	= &\frac{\lambda'-u}{F^2}\nabla_j\nabla_iF - \frac 1F\nabla_j\nabla_i(\l'-u)+\nabla_j\nabla_iu\nonumber\\
&+\frac 1F \(\nabla_jF\nabla_i(\frac{\lambda'-u}F)+\nabla_iF\nabla_j(\frac{\lambda'-u}F)\)+\(\frac{\lambda'-u}F-u\)\((h^2)_{ij}+g_{ij}\). \label{s3:h-evl}
	\end{align}
Substituting \eqref{2.1}, \eqref{2.2} and \eqref{s3:Sims} into \eqref{s3:h-evl}, we have
\begin{align}\label{s3:h-evl2}
	\frac{\partial}{\partial t} h_{ij}= & \frac{\lambda'-u}{F^2}\dot{F}^{kl}\nabla_k\nabla_lh_{ij}+\frac{\lambda'-u}{F^2}\ddot{F}^{kl,pq}\nabla_ih_{kl}\nabla_jh_{pq}+\frac 1F \(\nabla_jF\nabla_i(\frac{\lambda'-u}F)+\nabla_iF\nabla_j(\frac{\lambda'-u}F)\)\nonumber\\
&+\frac{\lambda'-u}{F^2}\dot{F}^{kl}((h^2)_{kl}+g_{kl})h_{ij}-\frac{\lambda'-u}{F^2}\dot{F}^{kl}h_{kl}((h^2)_{ij}+g_{ij})+\frac{\lambda'}F(h_{ij}-g_{ij})\nonumber\\
 &+\frac 1F\langle V,\nabla h_{ij}\rangle-\frac uF((h^2)_{ij}-h_{ij})+\langle V,\nabla h_{ij} \rangle +\l'  h_{ij} -u (h^2)_{ij}\nonumber\\
	& +\(\frac{\lambda'-u}F-u\)\((h^2)_{ij}+g_{ij}\).
	\end{align}
Let $S_{ij}=h_{ij}-g_{ij}$. Equations \eqref{s3:2-gij} and \eqref{s3:h-evl2} imply that
	\begin{align*}
	\frac{\partial}{\partial t} S_{ij}= & \frac{\lambda'-u}{F^2}\dot{F}^{kl}\nabla_k\nabla_lh_{ij}+\frac{\lambda'-u}{F^2}\ddot{F}^{kl,pq}\nabla_ih_{kl}\nabla_jh_{pq}+\frac 1F \(\nabla_jF\nabla_i(\frac{\lambda'-u}F)+\nabla_iF\nabla_j(\frac{\lambda'-u}F)\)\nonumber\\
&+\frac{\lambda'-u}{F^2}\dot{F}^{kl}((h^2)_{kl}+g_{kl})h_{ij}-\frac{\lambda'-u}{F^2}\dot{F}^{kl}h_{kl}((h^2)_{ij}+g_{ij})+\frac{\lambda'}F(h_{ij}-g_{ij})\nonumber\\
 &+\frac 1F\langle V,\nabla h_{ij}\rangle-\frac uF((h^2)_{ij}-h_{ij})+\langle V,\nabla h_{ij} \rangle +\l'  h_{ij} -u (h^2)_{ij}\nonumber\\
	& +\(\frac{\lambda'-u}F-u\)\((h^2)_{ij}-2h_{ij}+g_{ij}\)\nonumber\\
=&\frac{\lambda'-u}{F^2}\dot{F}^{kl}\nabla_k\nabla_lS_{ij}+\frac{\lambda'-u}{F^2}\ddot{F}^{kl,pq}\nabla_iS_{kl}\nabla_jS_{pq}+\frac 1F \(\nabla_jF\nabla_i(\frac{\lambda'-u}F)+\nabla_iF\nabla_j(\frac{\lambda'-u}F)\)\nonumber\\
 &+\frac {1+F}F\langle V,\nabla S_{ij}\rangle+\left(\frac{\lambda'-u}{F^2}(\dot{F}^{kl}(S^2)_{kl}+F)+\lambda'-2u\right)S_{ij}\nonumber\\
 &-\left(\frac{\lambda'-u}{F^2}\dot{F}^{kl}g_{kl}+\frac{u}F+2u\right)(S^2)_{ij}+(\lambda' -u) \left(1+\frac{\dot{F}^{kl}(S^2)_{kl}}{F^2}\right)g_{ij}.
	\end{align*}
	
	(iii) Since $F=F(\tilde{\kappa})$ is a smooth symmetric function of the shifted Weingarten matrix $S_i^j=S_{ik}g^{jk}$, by \eqref{s3:2-gij}, \eqref{s3:h-evl} and \eqref{2.1}, \eqref{2.2}, we have
	\begin{align*}
	\frac{\partial}{\partial t}F=&\frac{\partial F}{\partial S_{i}^l} \frac{\partial}{\partial t}S_{i}^{l}=\dot{F}^{ij}g_{jl}\frac{\partial}{\partial t}h_{i}^{l}\\
=& \dot{F}^{ij}\left(\frac{\partial}{\partial t}h_{ij}-h_{ik}g^{kl}\frac{\partial}{\partial t}g_{jl}\right)\\
= &\frac{\lambda'-u}{F^2}\dot{F}^{ij}\nabla_j\nabla_iF +\frac 2F\dot{F}^{ij}\nabla_jF\nabla_i(\frac{\lambda'-u}F)- \frac 1F\dot{F}^{ij}\nabla_j\nabla_i(\l'-u)\nonumber\\
&+\dot{F}^{ij}\nabla_j\nabla_iu-\(\frac{\lambda'-u}F-u\)\dot{F}^{ij}\((h^2)_{ij}-g_{ij}\)\nonumber\\
	=&\frac{\lambda'-u}{F^2}\dot{F}^{ij}\nabla_j\nabla_iF +\frac 2F\dot{F}^{ij}\nabla_jF\nabla_i(\frac{\lambda'-u}F)-\frac{\lambda'}F\dot{F}^{ij}((h^2)_{ij}-h_{ij})\nonumber\\
 &+\frac {1+F}F\dot{F}^{ij}\langle V,\nabla h_{ij}\rangle+\frac uF\dot{F}^{ij}(h_{ij}-g_{ij}) +\l'  \dot{F}^{ij}h_{ij} -u \dot{F}^{ij}g_{ij}\\
 =& \frac{\lambda'-u}{F^2}\dot{F}^{ij}\nabla_j\nabla_iF +\frac 2F\dot{F}^{ij}\nabla_jF\nabla_i(\frac{\lambda'-u}F)+\frac {1+F}F\langle V,\nabla F\rangle\nonumber\\
 &-\frac{\lambda'}F\dot{F}^{ij}(S^2)_{ij}+\l' F +(\lambda'-u) (\dot{F}^{ij}g_{ij}-1).
	\end{align*}

(iv) For the support function $u=\langle \lambda(r)\partial_r,\nu\rangle$, we first note that the unit outward normal $\nu$ evolves by
 \begin{equation*}
  \frac{\partial}{\partial t}\nu=-\nabla (\frac{\lambda'-u}F-u).
 \end{equation*}
Then using the property that $V=\lambda(r)\partial_r$ is a conformal Killing vector field, we have
\begin{align*}
  \frac{\partial}{\partial t}u = &\frac{\partial}{\partial t}\langle V,\nu\rangle\\
  =& \langle \lambda'(r) (\frac{\lambda'-u}F-u)\nu,\nu\rangle-\langle V, \nabla(\frac{\lambda'-u}F-u)\rangle\\
  =& \lambda'(\frac{\lambda'-u}F-u)-\langle V,\nabla(\frac{\lambda'-u}F-u)\rangle.
\end{align*}
By \eqref{2.2},
\begin{align*}
  \frac{\lambda'-u}{F^2}\dot{F}^{ij}\nabla_j\nabla_iu= & \frac{\lambda'-u}{F^2}\langle V,\nabla F \rangle+\frac{\lambda'(\lambda'-u)}{F^2}\dot{F}^{ij}(S_{ij}+g_{ij})\\
   &-\frac{u(\lambda'-u)}{F^2}\dot{F}^{ij}((S^2)_{ij}+2S_{ij}+g_{ij}) \\
   =&\frac{\lambda'-u}{F^2}\langle V,\nabla F \rangle-\frac{u(\lambda'-u)}{F^2}\dot{F}^{ij}(S^2)_{ij}+\frac{\lambda'-u}{F}(\lambda'-2u)\\
   &+\frac{(\lambda'-u)^2}{F^2}\dot{F}^{ij}g_{ij}.
\end{align*}
Combining the above two equations, we have
\begin{align*}
  \frac{\partial}{\partial t}u -\frac{\lambda'-u}{F^2}\dot{F}^{ij}\nabla_j\nabla_iu=& \lambda'(\frac{\lambda'-u}F-u)-\langle V,\nabla(\frac{\lambda'-u}F-u)\rangle-\frac{\lambda'-u}{F^2}\langle V,\nabla F \rangle \\
  &+\frac{u(\lambda'-u)}{F^2}\dot{F}^{ij}(S^2)_{ij}-\frac{\lambda'-u}{F}(\lambda'-2u)-\frac{(\lambda'-u)^2}{F^2}\dot{F}^{ij}g_{ij}\\
  =&\frac{1+F}F\langle V,\nabla u\rangle -\frac 1F\langle V,\nabla\lambda'\rangle+\frac{u(\lambda'-u)}{F^2}\dot{F}^{ij}(S^2)_{ij}\\
  &-\frac{(\lambda'-u)^2}{F^2}\dot{F}^{ij}g_{ij}-\lambda'u+2u\frac{\lambda'-u}{F}.
\end{align*}
	
(v) The time derivative of $\lambda'(r)$ satisfies
\begin{align}\label{s3:lam'-evl2}
  \frac{\partial}{\partial t}\lambda'(r)=&\langle V,\partial_tX\rangle =(\frac{\lambda'-u}F-u)u.
\end{align}
Combining the \eqref{s3:lam'-evl2} with \eqref{2.1} gives that
\begin{align*}
  \frac{\partial}{\partial t}\lambda'(r)-\frac{\lambda'-u}{F^2}\dot{F}^{ij}\nabla_i\nabla_j\lambda'(r)=&(\frac{\lambda'-u}F-u)u-\frac{\lambda'-u}{F^2}\dot{F}^{ij}(\lambda'g_{ij}-uS_{ij}-ug_{ij})\\
  =&2u\frac{\lambda'-u}F-u^2-\frac{(\lambda'-u)^2}{F^2}\dot{F}^{ij}g_{ij}.
\end{align*}
\end{proof}

\section{Preserving h-convexity}\label{sec:h-conv1}

We first recall the tensor maximum principle, which was first proved by Hamilton \cite{Ham1982} and was generalized by Andrews \cite{And07}.
\begin{thm}[Theorem 3.2 of \cite{And07}]\label{thm-2}
	Let $S_{ij}$ be a smooth time-varying symmetric tensor field on a compact manifold $M$ satisfying
	\begin{align*}
	\frac{\partial}{\partial t}S_{ij} = a^{k\ell} \nabla_k \nabla_\ell S_{ij}+u^k \nabla_k S_{ij}+N_{ij},
	\end{align*}
	where $a^{kl}$ and $u$ are smooth, $\nabla$ is a (possibly time-dependent) smooth symmetric connection, and $a^{kl}$ is positive definite everywhere. Suppose that
	\begin{align}\label{3.2}
	N_{ij} v^i v^j + \sup_{\L_k^p} 2 a^{k\ell}(2\L_k^p \nabla_\ell S_{ip}v^i-\L_k^p \L_\ell^q S_{pq}) \geq 0,
	\end{align}
whenever $S_{ij}\geq 0$ and $S_{ij}v^j=0$, where the supremum is taken over all $n\times n$ matrix $\Lambda_k^p$.  If $S_{ij}\geq 0$ everywhere on $M$ at $t=0$ and on $\partial M$ for $0\leq t\leq T$, then $S_{ij}\geq 0$ holds on $M$ for $0\leq t\leq T$.
\end{thm}
\begin{rem}
Note that in Andrews' paper, the conclusion is stated as ``if $S_{ij}$ is positive definite everywhere on $M$ at $t=0$ and on $\partial M$ for $0\leq t\leq T$, then $S_{ij}$ is positive definite for $0\leq t\leq T$. This in fact can be rephrased as the one stated in Theorem \ref{thm-2}. See the proof of Theorem 3.2 in Andrews' paper and Theorem 9.1 in Hamilton's paper \cite{Ham1982}.
\end{rem}

The idea in the proof of Theorem \ref{thm-2} is to consider the function on the tangent bundle of $M$:
\begin{equation*}
  Z(p,v)=S(p)(v,v),
\end{equation*}
where $p\in M$ and $v\in T_pM$. Let $p\in M$ be a point where $S(p)$ has a null vector $v$. To show that $Z\geq 0$ is preserved, it suffices to show that the time derivative of $Z$ at $(p,v)$ is nonnegative. Choose coordinates $\{x^i\}_{i=1}^n$ for $M$ near $p$ with connection coefficients for $\nabla$ vanish at $p$. Then the tangent bundle $TM$ near $p$ is described locally by $2n$ coordinates $\{x^1,\cdots, x^n, \dot{x}^1,\cdots, \dot{x}^n\}$. The null vector is $v=\partial_1$. The first order condition implied by minimality of $Z$ at $(p,v)$ is
\begin{equation*}
  0=\frac{\partial Z}{\partial x^i}=\frac{\partial Z}{\partial \dot{x}^i},\quad i=1,\cdots,n.
\end{equation*}
The second derivatives form a $2n\times 2n$ matrix which is non-negative at $(p,v)$. It follows that for any choice of matrix $\Lambda_k^p$, there holds
\begin{equation*}
  a^{k\ell}\left(\frac{\partial}{\partial x^k}-\L_k^p\frac{\partial}{\partial \dot{x}^p}\right)\left( \frac{\partial}{\partial x^\ell}-\L_\ell^q\frac{\partial}{\partial \dot{x}^q}\right)Z\geq 0.
\end{equation*}
By choosing the optimal $\L_k^p$ such that \eqref{3.2} holds, the time derivative $\partial Z/{\partial t}$ at $(p,v)$ is nonnegative. Note that the spatial minimal value of a function on $M$ is just a Lipschitz function in time. However, the argument in Lemma 3.5 of Hamilton's paper \cite{Ham1986} implies that the time derivative of $\min_{TM_t}Z(p,v)$ is not smaller than the value of $\partial Z/{\partial t}$ at the minimal point $(p,v)$, which implies that $Z\geq 0$ is preserved. Moreover, if strict inequality holds in \eqref{3.2}, we have that the minimum value of $Z$ over $TM_t$ is strictly increasing in time, and so $Z>0$ for later time.

Now we apply Theorem \ref{thm-2} to prove that the flow \eqref{s1:flow-BGL} preserves the h-convexity.
\begin{thm}\label{thm-h-convexity}
	If the initial hypersurface $M_0$ is h-convex, then along the flow \eqref{s1:flow-BGL}, the evolving hypersurface $M_t$ is strictly h-convex for $t>0$.
\end{thm}
\begin{proof}
Let $S_{ij}=h_{ij}-g_{ij}$. The h-convexity is equivalent to the positivity of the tensor $S_{ij}$. By \eqref{s3:1-gij} and \eqref{s3:1-hij}, we have
\begin{align*}
	\frac{\partial}{\partial t} S_{ij}=&\frac{\lambda'}{F^2}\dot{F}^{kl}\nabla_k\nabla_lh_{ij}+\frac{\lambda'}{F^2}\ddot{F}^{kl,pq}\nabla_ih_{kl}\nabla_jh_{pq}+\frac 1F\(\nabla_jF\nabla_i(\frac{\lambda'}F)+\nabla_iF\nabla_j(\frac{\lambda'}F)\)\nonumber\\
&+\langle V,\nabla h_{ij} \rangle+\(\frac{u}{F}+ \l'+\frac{\l'}{F^2}\dot{F}^{kl}((h^2)_{kl}+g_{kl})\) h_{ij} \nonumber\\
&-2u(h^2)_{ij}-\(\frac{\l'}{F}+u\)g_{ij}-2(\frac{\lambda'}F-u)h_{ij}\nonumber\\
=&\frac{\l'}{F^2}\dot{F}^{kl}\nabla_k\nabla_l S_{ij}+\frac{\l'}{F^2}\ddot{F}^{kl,pq}\nabla_ih_{kl}\nabla_jh_{pq}+\frac 1{F^2}\(\nabla_i\l'\nabla_jF+\nabla_j\l'\nabla_iF\)\\
	&-2\frac{\l'}{F^3}\nabla_iF\nabla_jF+\langle V,\nabla S_{ij} \rangle+\(\frac{u}{F}+\l'+\frac{\l'}{F^2}\dot{F}^{kl}((h^2)_{kl}+g_{kl})\) (S_{ij}+g_{ij})\\
&-2u((S^2)_{ij}+S_{ij})-2\frac{\l'}{F} (S_{ij}+g_{ij}) -\(\frac{\l'}{F}+u\)g_{ij}.
	\end{align*}
To apply the tensor maximum principle, we need to show that \eqref{3.2} whenever $S_{ij}\geq 0$ and $S_{ij}v^j=0$ (so that $v$ is a null vector of $S$). Let $(x_0,t_0)$ be the point where $S_{ij}$ has a null vector $v$. By continuity, we can assume that the principal curvatures are mutually distinct and in increasing order at $(x_0,t_0)$, that is $\k_1<\k_2<\cdots<\k_n$\footnote{This is possible since for any positive definite symmetric matrix $A$ with $A_{ij}\geq 0$ and $A_{ij}v^iv^j=0$ for some $v\neq 0$, there is a sequence of symmetric matrixes $\{A^{(k)}\}$ approaching $A$, satisfying $A^{(k)}_{ij}\geq 0$ and $A^{(k)}_{ij}v^iv^j=0$ and with each $A^{(k)}$ having distinct eigenvalues. Hence it suffices to prove the result in the case where all of $\kappa_i$ are distinct. }. The null vector condition $S_{ij}v^j=0$ implies that $v=e_1$ and $S_{11}=\k_1-1=0$ at $(x_0,t_0)$.	
	The terms involving $S_{ij}$ and $(S^2)_{ij}$ satisfy the null vector condition and can be ignored. Thus, it remains to show that
	\begin{align}\label{s4:Q1}
	Q_1:=&\frac{\l'}{F^2}\ddot{F}^{kl,rs}\nabla_1h_{kl}\nabla_1h_{rs}+\frac{2}{F^2}\nabla_1\l'\nabla_1F-2\frac{\l'}{F^3}|\nabla_1F|^2\nonumber\\
         &+\left( \frac{\l'}{F^2}(\dot{F}^{kl}((h^2)_{kl}+g_{kl})-2F)+(\l'-u)\frac{F-1}{F}\right)\nonumber\\
    	 &+2\frac{\l'}{F^2}  \sup_{\L} \dot{F}^{kl} (2\L_k^p \nabla_l S_{1p}-\L_k^p\L_l^q S_{pq}) \geq 0
	\end{align}
at $(x_0,t_0)$ for all matrix $(\L_k^p)$. By assumption, $S_{11}=0$ and $\nabla_k S_{11}=0$ at $(x_0,t_0)$. We have
	\begin{align*}
	\dot{F}^{kl} (2\L_k^p \nabla_l S_{1p}-\L_k^p\L_l^q S_{pq})=&\sum_{k=1}^{n}\sum_{p=2}^{n}\dot{F}^k(2\L_k^p \nabla_k S_{1p}-(\L_k^p)^2 S_{pp}) \nonumber \\
	=& \sum_{k=1}^{n} \sum_{p=2}^{n}\dot{F}^k \(\frac{(\nabla_k S_{1p})^2}{S_{pp}}-\(\L_k^p-\frac{\nabla_k S_{1p}}{S_{pp}}\)^2 S_{pp}\). 
	\end{align*}
Then the supremum of the last line in \eqref{s4:Q1} is obtained by choosing the matrix $\L_k^p=\frac{\nabla_k S_{1p}}{S_{pp}}$ for $p\geq 2$, $k\geq 1$ and $\L_k^1=0$ for all $k$. It follows that
 \begin{align}\label{s4:Q1-2}
	Q_1=&\frac{\l'}{F^2}\left(\ddot{F}^{kl,rs}\nabla_1h_{kl}\nabla_1h_{rs}+2\sum_{k=1}^{n}\sum_{p=2}^{n}\dot{F}^k \frac{(\nabla_k S_{1p})^2}{S_{pp}}\right)+\frac{2}{F^2}\nabla_1\l'\nabla_1F\nonumber\\
         &-2\frac{\l'}{F^3}|\nabla_1F|^2+\left( \frac{\l'}{F^2}(\dot{F}^{kl}((h^2)_{kl}+g_{kl})-2F)+(\l'-u)\frac{F-1}{F}\right).
	\end{align}
Moreover, since $F(\kappa)$ is inverse-concave\footnote{$F(\kappa)$ is inverse-concave, if its dual function $F_*(z)=F(\frac 1{z_1},\cdots,\frac 1{z_n})^{-1}$ is concave. As in \cite{And07}, the quotient $E_{m}(\kappa)/{E_{m-1}(\kappa)}$ is both concave and inverse concave.}, by a direct calculation as in \cite[\S 3]{AW18}, we have
	\begin{align}\label{s4:h-cov1}
        &\ddot{F}^{kl,rs}\nabla_1 h_{kl}\nabla_1  h_{rs}+2\sum_{k=1}^{n}\sum_{p=2}^{n}\dot{F}^k \frac{(\nabla_k S_{1p})^2}{S_{pp}} \nonumber \\
	\geq & \frac 2F|\nabla_1F|^2+2\sum_{k>1,l>1}\dot{F}^k\left(\frac{1}{\k_l-1}-\frac 1{\k_l}\right)(\nabla_1 h_{kl})^2 \nonumber \\
	\geq & \frac 2F|\nabla_1F|^2+2\sum_{k>1}\frac{\dot{F}^k}{\k_k(\k_k-1)}(\nabla_1 h_{kk})^2.
	\end{align}
	On the other hand, by the Cauchy-Schwarz inequality we have
	$$
	\sum_{k=2}^n\frac{\dot{F}^k}{\k_k(\k_k-1)}(\nabla_1 h_{kk})^2 \cdot \sum_{k=2}^n \dot{F}^k \k_k(\k_k-1) \geq \(\sum_{k=2}^n \dot{F}^k|\nabla_1 h_{kk}|\)^2 \geq |\nabla_1F|^2,
	$$
	where we used $\nabla_1 S_{11}=\nabla_1 h_{11}=0$ and \eqref{s2:0-2} that $\sum_{k=1}^n\dot{F}^k\k_k^2-F\geq F^2 -F>0$ since $1<\k_2<\cdots<\k_n$. Hence, we get
	\begin{align}\label{key-estimate-1}
	\sum_{k=2}^n\frac{\dot{F}^k}{\k_k(\k_k-1)}(\nabla_1 h_{kk})^2 \geq \frac{|\nabla_1F|^2}{\sum_{k=1}^n\dot{F}^k \k_k^2 -F}.
	\end{align}
	Since $\dot{F}^{kl}((h^2)_{kl}+g_{kl})-2F=\sum_{k=1}^{n}\dot{F}^k(\k_k-1)^2 \geq 0$, we estimate the last term of \eqref{s4:Q1-2} as follows
	\begin{align*}
	     &\frac{\l'}{F^2}(\dot{F}^{kl}((h^2)_{kl}+g_{kl})-2F)+(\l'-u)\frac{F-1}{F} \\
\geq & \frac{\l'-u}{F^2}(\sum_{k=1}^n\dot{F}^k\k_k^2+\sum_{k=1}^n\dot{F}^k-2F)+(\l'-u)\frac{F-1}{F}\\
	= &\frac{\l'-u}{F^2}\left(\sum_{k=1}^n\dot{F}^k\k_k^2+\sum_{k=1}^n\dot{F}^k-3F+F^2\right) \\
	\geq &\frac{\l'-u}{F^2}\left(\sum_{k=1}^n\dot{F}^k\k_k^2-F+(F-1)^2\right)\\
	\geq &\frac{\l'-u}{F^2}(\sum_{k=1}^n\dot{F}^k\k_k^2-F),
	\end{align*}
	where we used \eqref{s2:0-1} that $\sum_{k=1}^{n}\dot{F}^k\geq 1$ in the second inequality and the fact $\l'>u$. Therefore, we obtain
	\begin{align*}
	Q_1 \geq &\frac{2\nabla_1\l'\nabla_1F}{F^2}+\frac{2\l'}{F^2}\frac{|\nabla_1F|^2}{\sum_{k=1}^{n}\dot{F}^k \k_k^2-F}+\frac{\l'-u}{F^2}(\sum_{k=1}^n\dot{F}^k\k_k^2-F) \\
=&\frac{2\l'}{F^2(\sum_{k=1}^{n}\dot{F}^k \k_k^2-F)}\left(\nabla_1F+\frac {\nabla_1\lambda'}{2\lambda'}(\sum_{k=1}^n\dot{F}^k \k_k^2-F)\right)^2\\
&\quad -\frac {|\nabla_1\lambda'|^2}{2\lambda'F^2}(\sum_{k=1}^n\dot{F}^k \k_k^2-F)+\frac{\l'-u}{F^2}(\sum_{k=1}^n\dot{F}^k\k_k^2-F) \\
	    \geq &\frac{\sum_{k=1}^n\dot{F}^k\k_k^2-F}{2\lambda'F^2}\left(2\lambda'(\lambda'-u)-|\nabla_1\lambda'|^2\right).
	\end{align*}
Using \eqref{2.1}, we have the estimate $|\nabla_1\lambda'|^2=\langle V,e_1\rangle^2\leq \lambda^2-u^2$. Then
\begin{align*}
2\lambda'(\lambda'-u)-|\nabla_1\lambda'|^2
\geq & 2\lambda'(\lambda'-u)-(\lambda^2-u^2) \\
=&  (\l'-u)^2+\l'^2-\l^2 \geq 1.
\end{align*}
This implies that $Q_1$ is positive. Then the tensor maximum principle implies that $S_{ij}$ is positive for positive time $t>0$.
\end{proof}
\begin{rem}
By \cite[Lemmas 4 \& 5]{AMZ13}, if a smooth, symmetric and homogeneous of degree one function $F(\kappa)$ is concave and inverse-concave and is normalized such that $F(1,\cdots,1)=1$, then
\begin{equation*}
 \sum_{i=1}^{n}\dot{F}^i \geq 1,\qquad \sum_{i=1}^{n}\dot{F}^i \k_i^2 \geq F^2\quad\mathrm{ on}\quad \G_{+}=\{x\in \mathbb{R}^n,~x_i>0\}.
\end{equation*}
It's easy to see that the above proof of Theorem \ref{thm-h-convexity} works for such  functions $F(\kappa)$. Therefore, the h-convexity of the solution $M_t$ is preserved along the flow \eqref{s1:flow-BGL} with the quotient $E_m(\kappa)/{E_{m-1}(\kappa)}$ replaced by any concave and inverse concave homogeneous of degree one function $F(\kappa)$.
\end{rem}

\begin{proof}[Proof of Theorem \ref{main-thm-I}]
Firstly, the $C^0$ estimate of the flow \eqref{s1:flow-BGL} follows from the comparison principle directly.  Then the $C^1$ estimate follows from the h-convexity immediately (See e.g., the proof of \cite[Theorem 2.7.10]{Gerh06}, and the proof of Propositions \ref{C0-estimate} and \ref{pro-gradient-estimate}, which are similar results on $C^0, C^1$ estimates for the flow \eqref{s1:flow1}). With the $C^0, C^1$ estimates in hand, we can apply the maximum principle to the evolution equation \eqref{s3:1-hij}  to derive the upper bound on the principal curvatures (See also the proof of Proposition \ref{s8:prop-upper-bound} for the proof of curvature estimate for the flow \eqref{s1:flow1}). Then we have two sides positive bounds $1\leq \kappa_i\leq C$ on the principal curvatures along the flow \eqref{s1:flow-BGL}, which is equivalent to uniform $C^2$ estimate of the solution $M_t$ (as we already established the $C^0, C^1$ estimates). The smooth convergence to a geodesic sphere follows by estimating the evolution equation of the gradient of the radial function of $M_t$.
\end{proof}

\section{New geometric inequalities: I}\label{sec:Ineqn}

In this section, we will use the smooth convergence in Theorem \ref{main-thm-I} on the flow \eqref{s1:flow-BGL} to give the proof of Theorems \ref{thm-geom-inequality} and \ref{thm-weighted-AF-inequality}.

\subsection{Proof of Theorem \ref{thm-geom-inequality}}$\ $

We recall the following lemma which was proved in \cite[Lemma 4.5 \& Lemma 4.7]{GeWW14}.
\begin{lem}\label{s5:lem1}
	For any $\k\in \{x\in \mathbb R^n~|~ x_i>1 \}$, let
\begin{align*}
	\~L_k =\sum_{i=0}^{k}\binom ki (-1)^i E_{2k-2i}(\kappa), \qquad \~N_k =\sum_{i=0}^{k}\binom ki (-1)^i E_{2k-2i+1}(\kappa).
	\end{align*}
We have $\~L_k>0$, $\~N_k>0$ and
	\begin{align}\label{basic-Nk-Lk-inequality}
	E_{2k+1}(\kappa)\~L_k-E_{2k}(\kappa)\~N_k \geq 0.
	\end{align}
Equality holds if and only if $\k=c I$ for some constant $c>1$.
\end{lem}
We also recall the following variational formula (see \cite[Lemma 4.1]{GeWW14})
\begin{lem}
	Along any smooth variation
	\begin{align}\label{s5:flow}
	\frac{\partial}{\partial t} X(x,t)= \eta(x,t)\nu(x,t),
	\end{align}
	we have
	\begin{align}\label{variational-formula-L_k}
	\frac{d}{dt}\int_{M_t}\~L_k d\mu_t=(n-2k)\int_{M_t} \~N_k \eta d\mu_t.
	\end{align}
\end{lem}

By \eqref{variational-formula-L_k}, along the flow \eqref{s1:flow-BGL} we have
\begin{align*}
\frac{d}{dt}\int_{M_t}\~L_k  d\mu_t=&(n-2k)\int_{M_t} \~N_k \(\frac{\l'E_{m-1}}{E_m}-u\) d\mu_t \\
                              =&(n-2k)\int_{M_t} \l'\(\frac{E_{m-1}}{E_m}\~N_k-\~L_k\) d\mu_t,
\end{align*}
where we used the Minkowski identity \eqref{s2:MF} in the second equality. Since $1\leq m\leq 2k+1\leq n$, it follows from Lemma \ref{s5:lem1} and the Newton-MacLaurin inequality \eqref{s2:Newt-2} that
\begin{align}\label{crucial-inequality}
\frac{E_{m-1}}{E_{m}}\~N_k -\~L_k\leq \frac{E_{2k}}{E_{2k+1}}\~N_k -\~L_k \leq 0.
\end{align}
Then
 \begin{equation*}
  \frac{d}{dt}\int_{M_t}\~L_k d\mu_t\leq 0,
 \end{equation*}
and $\int_{M_t}\~L_kd\mu_t$ is monotone non-increasing.

On the other hand, the quermassintegral $W_m(\Omega_t)$ remains to be a constant along the flow \eqref{s1:flow-BGL}. By Theorem \ref{main-thm-I}, $M_t$ smoothly converges to a geodesic sphere of radius $r_{\infty}$ satisfying that
\begin{align*}
W_m(\Omega_0)=W_m(B_{r_{\infty}})=f_m(r_{\infty}).
\end{align*}
Since $L_k=\binom n{2k}(2k)!\tilde{L}_{k}$, the monotonicity of $\int_{M}\~L_k$ yields
\begin{align*}
\int_{M_0}L_k d\mu_0 \geq &\int_{\partial B_{r_{\infty}}} L_k d\mu_{\infty}=\binom n{2k}(2k)! \omega_{n}\sinh^{n-2k}r_{\infty} \\
= &g_k(r_{\infty})=g_k \circ f_m^{-1} (W_m(\Omega_0)).
\end{align*}
This gives the inequality \eqref{GWW-inequality} for smooth h-convex domains.

It is easy to check that equality holds in \eqref{GWW-inequality} for geodesic spheres. To show the converse, suppose that a smooth h-convex domain $\Omega$ achieves the equality of \eqref{GWW-inequality}. Then along the flow \eqref{s1:flow-BGL}, the integral $\int_{M_t}\~L_k d\mu_t$ remains to be a constant. It follows from \eqref{crucial-inequality} that
\begin{align*}
\frac{E_{m-1}}{E_m}\~N_k-\~L_k \equiv 0, \quad \text{ on $M_t$}.
\end{align*}
Since for each $t>0$, the hypersurface $M_t=\partial\Omega_t$ is strictly h-convex, by the equality characterization of \eqref{basic-Nk-Lk-inequality} we have $M_t$ is totally umbilical and hence a geodesic sphere. Since the initial hypersurface $M_0$ is smoothly approximated by a family of geodesic sphere, it is also a geodesic sphere. This completes the proof of Theorem \ref{thm-geom-inequality}.

\subsection{Proof of Theorem \ref{thm-weighted-AF-inequality}}$\ $

Now we give the proof of Theorem \ref{thm-weighted-AF-inequality}. The key point is to establish the monotonicity of $\int_{M_t} \lambda'(r) E_k$ along the flow \eqref{s1:flow-BGL}, where $\lambda'(r)=\cosh r$. We first calculate the evolution of $\int_{M_t}\l' E_kd\mu_t$ along a general flow equation \eqref{s5:flow}. Firstly, the area element evolves as
\begin{equation}\label{s5:dmu}
  \frac{\partial}{\partial t}d\mu_t=H\eta d\mu_t
\end{equation}
along the flow \eqref{s5:flow}, where $H=nE_1(\kappa)$ is the mean curvature of $M_t$. The function $\lambda'(r)|_{M_t}$ evolves by
\begin{align}\label{s5:lamb-evl}
   \frac{\partial}{\partial t}\lambda'(r) =& \langle \bar{\nabla}\lambda'(r), \frac{\partial}{\partial t}X\rangle\nonumber\\
  =& \langle \lambda(r)\partial_r,\eta\nu\rangle\nonumber\\
  =&\eta u.
\end{align}
Moreover, by the similar calculation as in (iii) of Lemma \ref{s3:lem2},
\begin{align}\label{s5:2Ek}
  \frac{\partial}{\partial t}E_k(\kappa) =& \dot{E}_k^{ij}\left(-\nabla_i\nabla^j \eta-\eta((h^2)_{i}^j-\delta_{i}^j)\right),
\end{align}
where $\dot{E}_k^{ij}$ denotes the derivative of $E_k$ with respect to the Weingarten matrix. Combining \eqref{s5:dmu} -- \eqref{s5:2Ek}, we have
	\begin{align*}
	\frac{d}{dt}\int_{M_t}\l' E_k(\kappa) d\mu_t=&\int_{M_t} u E_k \eta  d\mu_t+\int_{M_t}\l' E_k H\eta  d\mu_t\\
& +\int_{M_t}\l' \dot{E}_k^{ij} \left(-\nabla_i\nabla^j \eta-\eta((h^2)_{i}^j-\delta_{i}^j)\right)  d\mu_t.
	\end{align*}
Since $\dot{E}_k^{ij} $ is divergence-free, we have
	\begin{align*}
	\frac{d}{dt}\int_{M_t}\l' E_k(\kappa) d\mu_t=&\int_{M_t} \left(u E_k+n\l'E_k E_1-\l'\dot{E}_k^{ij}((h^2)_i^j-\d_i^j)\right) \eta  d\mu_t-\int_{M_t}\nabla_i \nabla^j(\l') \dot{E}_k^{ij} \eta  d\mu_t \\
	=&\int_{M_t} \left(u E_k+n\l'E_k E_1-\l'\dot{E}_k^{ij} ((h^2)_i^j-\d_i^j)\right) \eta  d\mu_t -\int_{M_t} \dot{E}_k^{ij}(\l'\d_i^j-u h_i^j)\eta  d\mu_t\\
	=&\int_{M_t} \left(u E_k+n\l'E_k E_1-\l'\dot{E}_k^{ij} (h^2)_i^j +u\dot{E}_k^{ij} h_i^j\right) \eta  d\mu_t\\
	=&\int_{M_t} \biggl( u E_k+n\l'E_k E_1-\l'(nE_1E_k-(n-k)E_{k+1})+kuE_k \biggr) \eta  d\mu_t\\
	=&\int_{M_t} \biggl( (k+1)u E_k+(n-k)\l'E_{k+1}\biggr) \eta  d\mu_t,
	\end{align*}
where we used \eqref{newton-formula-1}, \eqref{newton-formula-3} and \eqref{2.1}. Note that when $k=n$, we set $E_{n+1}(\kappa)=0$.

For each integer $m=1,\cdots,n$,  if we choose
\begin{equation*}
  \eta=\frac{\l'E_{m-1}}{E_m}-u,\quad 1\leq m\leq k,
\end{equation*}
then
	\begin{align*}
	\frac{d}{dt}\int_{M_t}\l' E_k(\kappa)d\mu_t=&\int_{M_t} \left((k+1)u E_k+(n-k)\l'E_{k+1}\right)\(\frac{\l'E_{m-1}}{E_m}-u\) d\mu_t\\
	=&\int_{M_t} \biggl((k+1)u \(\l'\frac{E_{k}E_{m-1}}{E_m}-uE_k\)\\
&+(n-k)\l'\(\l'\frac{E_{k+1}E_{m-1}}{E_m}-u E_{k+1}\)\biggr) d\mu_t\\
	\leq &\int_{M_t}  \biggl((k+1)u \(\l'E_{k-1}-uE_k\)+(n-k)\l'\(\l'E_k-u E_{k+1}\)\biggr) d\mu_t,
	\end{align*}
where we used the Newton-MacLaurin inequality \eqref{s2:Newt-2}. Using \eqref{2.1} and \eqref{newton-formula-1}--\eqref{newton-formula-2}, we have
\begin{align}\label{s5:2-1}
  \dot{E}_k^{ij} \nabla_i\nabla^j\lambda'(r)= & \dot{E}_k^{ij}(\lambda'\delta_i^j-uh_i^j)=k(\lambda'E_{k-1}-uE_k).
\end{align}
Then
	\begin{align}\label{s5:2dt}
	\frac{d}{dt}\int_{M_t}\l' E_k(\kappa)d\mu_t \leq & \int_{M_t}  \left(\frac{k+1}{k}u \dot{E}_k^{ij} \nabla_i\nabla^j\lambda' +\frac{n-k}{k+1}\l'\dot{E}_{k+1}^{ij}\nabla_i\nabla^j\lambda' \right)  d\mu_t \nonumber\\
	=& -\int_{M_t} \left(\frac{k+1}{k} \dot{E}_k^{ij} \nabla^j \lambda' \nabla_i u +\frac{n-k}{k+1}\dot{E}_{k+1}^{ij}\nabla^j\lambda' \nabla_i \l'\right) d\mu_t \nonumber\\
	=& -\int_{M_t} \left(\frac{k+1}{k} \dot{E}_k^{ij} h_i^l \nabla^j \lambda' \nabla_l \lambda' +\frac{n-k}{k+1}\dot{E}_{k+1}^{ij}\nabla^j\lambda' \nabla_i \lambda' \right)  d\mu_t,
	\end{align}
where we used \eqref{2.1}, \eqref{2.2} that $\nabla_i u=h_i^l\nabla_l\lambda'$. Since $M_t$ is h-convex, all the matrices $\dot{E}_k^{ij}$, $\dot{E}_{k+1}^{ij}$ and $h_i^l$ are all positive definite. This implies that
\begin{equation*}
  \frac{d}{dt}\int_{M_t}\l' E_k(\kappa)d\mu_t\leq 0.
\end{equation*}
Since the flow \eqref{s1:flow-BGL} preserves the $m$th quermassintegral $W_m(\Omega_t)$ and converges to a geodesic sphere, we obtain the inequality \eqref{weighted-AF-inequality} for $1\leq m \leq k$. Inequality \eqref{weighted-AF-inequality} with $m=0$ follows immediately from the quermassintegral inequality \eqref{WX-inequality} and the fact that $h_k$ is strictly increasing in $r$. If equality holds in \eqref{weighted-AF-inequality}, we obtain from \eqref{s5:2dt} that $\nabla\lambda'(r)=0$ everywhere on $M_t$ for all $t$, and the initial hypersurface is a geodesic sphere centered at the origin.  This completes the proof of Theorem \ref{thm-weighted-AF-inequality}.

\section{New locally constrained curvature flow}\label{sec:new flow}
In the rest of this paper, we study the new locally constrained curvature flow \eqref{s1:flow1}. We write
\begin{equation*}
 F(\tilde{\kappa})= \frac{E_m(\tilde{\kappa})}{E_{m-1}(\tilde{\kappa})},
\end{equation*}
where $\tilde{\kappa}=(\kappa_1-1,\cdots,\kappa_n-1)$ are shifted principal curvatures, defined as the eigenvalues of the shifted Weingarten matrix $\mathcal{W}-I$. Then the flow \eqref{s1:flow1} can be written as
\begin{equation}\label{s6:flow1}
 \frac{\partial}{\partial t}X=\left(\frac{\lambda'-u}{F(\tilde{\kappa})}-u\right)\nu.
\end{equation}
We assume that the initial hypersurface $M_0$ is h-convex. Then it is star-shaped with respect to a point inside the domain $\Omega_0$ enclosed by $M_0$.  As discussed in \S \ref{sec:2-3}, we can equivalently write \eqref{s6:flow1} as a scalar parabolic PDE on $\mathbb{S}^n$ for the radial graph function $r(\cdot,t)$ and also for $\varphi(\cdot,t)$. By \eqref{s2:gamma-evl}, the function $\varphi\in C^{\infty}(\mathbb{S}^n\times[0,T))$ satisfies
\begin{equation}\label{s6:varphi}
  \frac{\partial}{\partial t}\varphi=\left(\frac{\lambda'-u}{F(\tilde{\kappa})}-u\right)\frac v{\lambda},
\end{equation}
where $v=\sqrt{1+|D\varphi|^2}$ is the function defined in \eqref{s2:v-def}.

We first prove the $C^0$-estimate of of solution to \eqref{s1:flow1}. This is equivalent to the $C^0$-estimate of $\varphi$ defined in \S \ref{sec:2-3}.
\begin{prop}\label{pro-C0-estimate}
	Let $\varphi\in C^{\infty}(\mathbb{S}^n\times[0,T))$ be a smooth solution to the initial value problem \eqref{s6:varphi}. Then
	\begin{align}\label{C0-estimate}
	\min_{\t \in\mathbb S^n} \varphi_0(\t)\leq \g(\t,t) \leq  \max_{\t\in \mathbb S^n} \varphi_0(\t),  \quad \forall (\t,t)\in \mathbb{S}^n \times [0,T).
	\end{align}
\end{prop}
\begin{proof}
	The proof is by the standard maximum principle. Here we only prove the upper bound, and the lower bound can be proved in a similar manner. At the spatial maximum point of $\g$, we have
	\begin{align*}
	D \g=0 ,\qquad  D^2 \g \leq 0.
	\end{align*}
	Then we have
\begin{equation*}
  v=1,\quad  u=\lambda, \quad h_i^j \geq \frac{\l'}{\l}\d_i^j
\end{equation*}
 at the spatial maximum point of $\varphi$. This implies that
\begin{align*}
  (\frac{\lambda'-u}F-u)\frac v{\lambda} = & ~\frac 1{\lambda}(\frac{\lambda'-\lambda}{F(\tilde{\kappa})}-\lambda ) \\
  \leq  &~ \frac 1{\lambda}(\frac{\lambda'-\lambda}{F((\frac{\lambda'}{\lambda}-1)I)}-\lambda )\\
  =& ~0,
\end{align*}
where we used the $1$-homogeneity of $F$ and the fact that $F$ is strictly increasing in each argument. Therefore,
\begin{align*}
  \frac d{dt}\max_{\mathbb{S}^n}\varphi(\cdot,t) \leq & ~0
\end{align*}
and hence $\max_{\t\in \mathbb S^n}\g(\t,t)\leq \max_{\t\in \mathbb S^n}\varphi_0(\t)$.
\end{proof}
\begin{rem}
	Though $\max\varphi$ is just a Lipschitz function, we can still get an estimate on its time derivative by using an observation due to Hamilton \cite[Lemma 3.5]{Ham1986}.
\end{rem}

\section{Curvature estimate}\label{sec:curv}
In this section, we first prove that the function $F(\tilde{\kappa})$ is bounded from above and below by positive constants. Then we show that the h-convexity of the initial hypersurface is preserved along the flow \eqref{s1:flow1}. Finally, using the bounds on $F(\tilde{\kappa})$ and applying maximum principle to the evolution equation of the shifted Weingarten matrix, we prove the uniform upper bound on the shifted principal curvatures.

\subsection{Estimate on $F$}$\ $

We first estimate the lower bound on $F$. By assumption, the shifted principal curvatures $\tilde{\kappa}\in \Gamma_m^+$ on the initial hypersurface. This combined with the compactness of the hypersurface gives a positive lower bound on $F(\tilde{\kappa})$ on $M_0$. Since the zero order terms of the evolution equation \eqref{evol-F} of $F(\tilde{\kappa})$ have no desired sign, we can not obtain the positive lower bound on $F$ on $M_t$ for positive time using the maximum principle directly.

\begin{lem}\label{s4:F-lbd}
There exists a constant $C>0$ depending only on $M_0$ such that $F\geq C>0$ along the flow \eqref{s1:flow1}.
\end{lem}
\proof
We consider the function
\begin{equation*}
  \psi=\frac{\lambda'-u}F+ \lambda',
\end{equation*}
which is bounded from above initially.  We need to derive the evolution equation of $\omega$ along the flow \eqref{s1:flow1}. Combining \eqref{s3:u-evl} and \eqref{s3:lambd'-evl}, we have
\begin{align}\label{s3:lamb-u-evl}
  \frac{\partial}{\partial t}(\lambda'-u) =&\frac{\lambda'-u}{F^2}\dot{F}^{ij}\nabla_j\nabla_i(\lambda'-u)+\frac{1+F}F\langle V,\nabla (\lambda'-u)\rangle -\langle V,\nabla\lambda'\rangle\nonumber\\
  &-\frac{u(\lambda'-u)}{F^2}\dot{F}^{ij}(S^2)_{ij}+(\lambda'-u)u.
\end{align}
Then \eqref{evol-F}, \eqref{s3:lambd'-evl}  and \eqref{s3:lamb-u-evl} imply that
\begin{align}\label{s7:omeg-evl}
   &\frac{\partial}{\partial t}\psi-\frac{\lambda'-u}{F^2}\dot{F}^{ij}\nabla_i\nabla_j\psi-\frac{1+F}F\langle V,\nabla \psi\rangle \nonumber\\
  =&-\frac{2+F}F{\langle V,\nabla\lambda'\rangle}+\frac{(\lambda'-u)^2}{F^3}\dot{F}^{ij}(S^2)_{ij}\nonumber\\
   &-\frac{(\lambda'-u)^2}F -\frac{(\lambda'-u)^2}{F^2} (2\dot{F}^{ij}g_{ij}-1)-u^2+2 u\frac{\lambda'-u}F.
\end{align}
Note that $\langle V,\nabla\lambda'\rangle=\lambda^2-u^2$. Applying \eqref{s2:0-1} and \eqref{s2:0-2} to \eqref{s7:omeg-evl}, we have
\begin{align}\label{s7:omeg-evl2}
  & \frac{\partial}{\partial t}\psi-\frac{\lambda'-u}{F^2}\dot{F}^{ij}\nabla_i\nabla_j\psi-\frac{1+F}F\langle V,\nabla \psi\rangle \nonumber\\
 \leq   &~-\frac{2+F}F(\lambda^2-u^2)+(n-m)\frac{(\lambda'-u)^2}F\nonumber\\
  & -\frac{(\lambda'-u)^2}{F^2} -u^2+2 u\frac{\lambda'-u}F\nonumber\\
   \leq & -\frac{(\lambda'-u)^2}{F^2}+\frac{\lambda'-u}F\left(2u+(n-m)(\lambda'-u)\right)-\lambda^2\nonumber\\
   =& -\psi^2+\big(2\lambda'+2u+(n-m)(\lambda'-u)\big)\psi\nonumber\\
   &-(\lambda')^2-\lambda^2-\lambda'\left(2u+(n-m)(\lambda'-u)\right).
\end{align}
By the $C^0$ estimate, both the last line and the coefficient of $\psi$ on the right hand side of \eqref{s7:omeg-evl2} are uniformly bounded. Then the maximum principle implies that $\psi$ is bounded from above. This together with the $C^0$ estimate implies the uniform positive lower bound on $F$.
\endproof

We next prove the upper bound on $F$.
\begin{lem}\label{s7:lem1}
There exists a constant $C>0$ depending only on $M_0$ such that $F\leq C$ along the flow \eqref{s1:flow1}.
\end{lem}
\proof
We consider the function
\begin{equation*}
  \xi=\ln F+2\lambda'-u.
\end{equation*}
The equation \eqref{evol-F} implies that
\begin{align}\label{s7:lnF}
&	\frac{\partial}{\partial t}\ln F-\frac{\lambda'-u}{F^2}\dot{F}^{ij}\nabla_j\nabla_i\ln F-\frac {1+F}F\langle V,\nabla \ln F\rangle\nonumber\\
=& -\frac {\lambda'-u} {F^4}\dot{F}^{ij}\nabla_jF\nabla_iF+\frac {2} {F^3}\dot{F}^{ij}\nabla_jF\nabla_i(\lambda'-u)-\frac{\lambda'}{F^2}\dot{F}^{ij}(S^2)_{ij}\nonumber\\
 &+\l' +\frac{(\lambda'-u)}F (\dot{F}^{ij}g_{ij}-1).
	\end{align}
Combining \eqref{s3:u-evl} and \eqref{s3:lambd'-evl}, we have
\begin{align}\label{s7:2l-u}
  &\frac{\partial}{\partial t}(2\lambda'-u) -\frac{\lambda'-u}{F^2}\dot{F}^{ij}\nabla_j\nabla_i(2\lambda'-u)-\frac{1+F}F\langle V,\nabla (2\lambda'-u)\rangle \nonumber\\
  =& -(\frac 1F+2)\langle V,\nabla\lambda'\rangle+4u\frac{\lambda'-u}F-2u^2-\frac{2(\lambda'-u)^2}{F^2}\dot{F}^{ij}g_{ij}\nonumber\\
  &-\frac{u(\lambda'-u)}{F^2}\dot{F}^{ij}(S^2)_{ij}+\frac{(\lambda'-u)^2}{F^2}\dot{F}^{ij}g_{ij}+\lambda'u-2u\frac{\lambda'-u}{F}\nonumber\\
  =& -(\frac 1F+2)(\lambda^2-u^2) +2u\frac{\lambda'-u}F-2u^2-\frac{(\lambda'-u)^2}{F^2}\dot{F}^{ij}g_{ij}\nonumber\\
  &-\frac{u(\lambda'-u)}{F^2}\dot{F}^{ij}(S^2)_{ij}+\lambda'u,
\end{align}
where we used $\langle V,\nabla\lambda'\rangle=\lambda^2-u^2$. By \eqref{s7:lnF} and \eqref{s7:2l-u} we have
\begin{align}\label{s7:omeg}
   & \frac{\partial}{\partial t}\xi -\frac{\lambda'-u}{F^2}\dot{F}^{ij}\nabla_j\nabla_i\xi-\frac{1+F}F\langle V,\nabla\xi\rangle\nonumber\\
  = & -\frac {\lambda'-u} {F^4}\dot{F}^{ij}\nabla_jF\nabla_iF+\frac {2} {F^3}\dot{F}^{ij}\nabla_jF\nabla_i(\lambda'-u)\nonumber\\
  &-\frac 1{F^2}\left(\lambda'+u(\lambda'-u)\right)\dot{F}^{ij}(S^2)_{ij}\nonumber\\
 &+\l' +\frac{(\lambda'-u)}F (\dot{F}^{ij}g_{ij}-1) -(\frac 1F+2)(\lambda^2-u^2)\nonumber\\
  &+2u\frac{\lambda'-u}F-2u^2-\frac{(\lambda'-u)^2}{F^2}\dot{F}^{ij}g_{ij}+\lambda'u.
\end{align}
At the spatial maximum point of $\xi$, we have $\nabla_i(\lambda'-u)=-\nabla_i\ln F-\nabla_i\lambda'$. Using \eqref{s2:0-1} and \eqref{s2:0-2}, the equation \eqref{s7:omeg} implies that
\begin{align}\label{s7:omeg2}
  &\frac{\partial}{\partial t} \xi -\frac{\lambda'-u}{F^2}\dot{F}^{ij}\nabla_j\nabla_i \xi-\frac{1+F}F\langle V,\nabla \xi\rangle\nonumber\\
  \leq  & -\frac {\lambda'-u} {F^4}\dot{F}^{ij}\nabla_jF\nabla_iF-\frac {2} {F^4}\dot{F}^{ij}\nabla_iF\nabla_jF-\frac {2} {F^3}\dot{F}^{ij}\nabla_jF\nabla_i\lambda'\nonumber\\
  &+\frac 1F\left( (m-1)(\lambda'-u)-(\lambda^2-u^2)+2u(\lambda'-u)\right) \nonumber\\
  &-2\lambda^2+u^2-\frac{(\lambda'-u)^2}{F^2}\nonumber\\
  \leq &-\lambda^2+\frac 1{2F^2}\dot{F}^{ij}\nabla_i\lambda'\nabla_j\lambda'-\frac{(\lambda'-u)^2}{F^2}\nonumber\\
  &+\frac 1F\left( (m-1)(\lambda'-u)-(\lambda^2-u^2)+2u(\lambda'-u)\right) \nonumber\\
    \leq &-\lambda^2+\frac {\lambda^2-u^2}{2F^2}\dot{F}^{ij}g_{ij}-\frac{(\lambda'-u)^2}{F^2}\nonumber\\
  &+\frac 1F\left( (m-1)(\lambda'-u)-(\lambda^2-u^2)+2u(\lambda'-u)\right).
\end{align}
By the $C^0$ estimate, $\lambda\geq C>0$. When $\xi$ is sufficiently large, we have $F$ is sufficiently large as well. In this case, the right hand side of \eqref{s7:omeg2} is non-positive. The maximum principle implies that $\xi\leq C$ along the flow \eqref{s1:flow1}. This in turn gives the upper bound on $F$.
\endproof

\subsection{Preserving h-convexity}\label{sec:h-conv2}$\ $

In this subsection, we apply the tensor maximum principle in Theorem \ref{thm-2} to prove that the flow \eqref{s1:flow1} preserves the h-convexity.
\begin{thm}\label{thm-h-convexity-2}
If the initial hypersurface $M_0$ is h-convex, then along the flow \eqref{s1:flow1} the evolving hypersurface $M_t$  is strictly h-convex for positive time.
\end{thm}
\begin{proof}
The h-convexity is equivalent to the positivity of the shifted second fundamental form $S_{ij}=h_{ij}-g_{ij}$. By \eqref{s3:2-hij} we have
	\begin{align*}
	\frac{\partial}{\partial t} S_{ij}
=&\frac{\lambda'-u}{F^2}\dot{F}^{kl}\nabla_k\nabla_lS_{ij}+\frac{\lambda'-u}{F^2}\ddot{F}^{kl,pq}\nabla_iS_{kl}\nabla_jS_{pq}\nonumber\\
 &+\frac 1F\(\nabla_jF\nabla_i(\frac{\lambda'-u}F)+\nabla_i F\nabla_j(\frac{\lambda'-u}F)\)\nonumber\\
&+\frac {1+F}F\langle V,\nabla S_{ij}\rangle+\left(\frac{\lambda'-u}{F^2}(\dot{F}^{kl}(S^2)_{kl}+F)+\lambda'-2u\right)S_{ij}\nonumber\\
&-\left(\frac{\lambda'-u}{F^2}\dot{F}^{kl}g_{kl}+\frac{u}F+2u\right)(S^2)_{ij}+(\lambda' -u) \left(1+\frac{\dot{F}^{kl}(S^2)_{kl}}{F^2}\right)g_{ij}.
	\end{align*}
To apply the tensor maximum principle, we need to show that \eqref{3.2} holds whenever ${S}_{ij}\geq 0$ and ${S}_{ij}v^j=0$. By assumption $({S}_{ij})\geq 0$ at the initial time. Assume that $({S}_{ij})\geq 0$ on $M\times [0,t_0]$ and there exists a point $x_0\in M_{t_0}$ and a direction $v\in T_{x_0}M_{t_0}$ such that ${S}_{ij}v^j=0$ at $(x_0,t_0)$. By continuity we can assume that the principal curvatures are mutually distinct at $(x_0,t_0)$ and in increasing order $\kappa_1<\kappa_2<\cdots<\kappa_n$. The null vector condition ${S}_{ij}v^j=0$ implies that $v=e_1$ and ${S}_{11}=\k_1-1=0$ at $(x_0,t_0)$. The terms involving ${S}_{ij}$ and $({S}^2)_{ij}$ satisfy the null vector condition and can be ignored.

It suffices to show that
	\begin{align}\label{s7:Q1}
	Q_1:=&\frac{\lambda'-u}{F^2}\ddot{F}^{kl,pq}\nabla_1S_{kl}\nabla_1S_{pq}+\frac 2F\nabla_1F\nabla_1(\frac{\lambda'-u}F)\nonumber\\
 &+(\lambda' -u) \left(1+\frac{\dot{F}^{kl}(S^2)_{kl}}{F^2}\right)\nonumber\\
         &+2\frac{\lambda'-u}{F^2}\sup_{\L_k^p} \dot{F}^{kl} (2\L_k^p \nabla_l {S}_{1p}-\L_k^p\L_l^q {S}_{pq}) \geq 0.
	\end{align}
By assumption,  ${S}_{11}=0$ and $\nabla_k {S}_{11}=0$ at $(x_0,t_0)$. As in the proof of Theorem \ref{thm-h-convexity},  the supremum in the last line of \eqref{s7:Q1} is obtained by choosing $\L_k^p=\frac{\nabla_k {S}_{1p}}{{S}_{pp}}$ for $p\geq 2, k\geq 1$ and $\L_k^1=0$ for all $k$, and we have
\begin{align}\label{s7:Q1-2}
	Q_1=&\frac{\lambda'-u}{F^2}\left(\ddot{F}^{kl,pq}\nabla_1S_{kl}\nabla_1S_{pq}+2\sum_{k=1}^{n}\sum_{p=2}^{n}\dot{F}^k \frac{(\nabla_k {S}_{1p})^2}{{S}_{pp}}\right)\nonumber\\
 &+\frac 2F\nabla_1F\nabla_1(\frac{\lambda'-u}F)+(\lambda' -u) \left(1+\frac{\dot{F}^{kl}(S^2)_{kl}}{F^2}\right).
\end{align}
Using the inverse-concavity of $F$ and a similar calculation as in \cite[\S 3]{AW18}, we have
	\begin{align}\label{s7:h-cov1}
        \ddot{F}^{kl,pq}\nabla_1 S_{kl}\nabla_1  S_{pq}+2\sum_{k=1}^{n}\sum_{p=2}^{n}\dot{F}^k \frac{(\nabla_k {S}_{1p})^2}{{S}_{pp}}\geq & \frac 2F|\nabla_1F|^2.
	\end{align}
Note that the estimate \eqref{s7:h-cov1} is different with \eqref{s4:h-cov1}, as here $F(\tilde{\kappa})$ is a function of the shifted principal curvatures $\tilde{\kappa}$. For the second term of \eqref{s7:Q1-2}, using \eqref{2.1} -- \eqref{2.2} we have
\begin{align*}
  \frac 2F\nabla_1F\nabla_1(\frac{\lambda'-u}F)= & \frac 2{F^2}\nabla_1F \nabla_1(\lambda'-u)-\frac{2(\lambda'-u)}{F^3}|\nabla_1F|^2\\
  = & \frac 2{F^2}\nabla_1F \langle V,e_1\rangle (1-\kappa_1)-\frac{2(\lambda'-u)}{F^3}|\nabla_1F|^2\\
  = &-\frac{2(\lambda'-u)}{F^3}|\nabla_1F|^2.
\end{align*}
Then the quantity $Q_1$ satisfies
\begin{align*}
	Q_1\geq &(\lambda' -u) \left(1+\frac{\dot{F}^{kl}(S^2)_{kl}}{F^2}\right) \geq  2(\lambda' -u)>0,
	\end{align*}
where we used the estimate \eqref{s2:0-2}. Then the minimum of the smallest shifted principal curvature over $M_t$ is strictly increasing whenever it approaches zero, and we arrive at the conclusion of Theorem \ref{thm-h-convexity-2}.
\end{proof}

\subsection{Curvature estimate}$\ $

We now prove the upper bound of the shifted principal curvatures.
\begin{prop}\label{s8:prop-upper-bound}
Along the flow \eqref{s1:flow1}, the shifted principal curvatures $M_t$ satisfy
\begin{equation*}
  \tilde{\kappa}_i\leq C,\qquad i=1,\cdots,n.
\end{equation*}
for some positive constant $C$.
\end{prop}
\proof
Define a function
\begin{equation*}
  \zeta(x,t)=\sup\{S_{ij}\eta^i\eta^j:~g_{ij}\eta^i\eta^j=1\},
\end{equation*}
which is the largest shifted principal curvature of $M_t$ at $(x,t)$. That is, $\zeta(x,t)$ is the largest eigenvalue of shifted Weingarten matrix $(S_i^j)$ at $(x,t)$. To estimate the upper bound on $\zeta$, we first need to calculate the evolution equation of $\zeta$. Basically this is just the evolution of $S_n^n$ with $e_n$ being the direction corresponding to the largest shifted principal curvature.

We calculate the evolution of $\zeta$ in details. For any time $t_0\in [0,T)$, we assume that the maximum of $\zeta(\cdot,t_0)$ is achieved at the point $x_0$ in the direction $\eta=e_n$, where we choose an orthonormal basis $\{e_1,\cdots,e_n\}$ at $(x_0,t_0)$.  We may extend $\{e_i\}_{i=1}^n$ to a space-time neighbourhood of $(x_0,t_0)$ by first parallel translating in space along radial geodesics emanating from $x_0$ with respect to $g(t_0)$, and then extending in the time direction by solving
 \begin{equation}\label{s7.3pf1}
   \frac {\partial}{\partial t}e_i(t)=-\left(\frac{\lambda'(r)-u}{F(\tilde{\kappa})}-u\right) \mathcal{W} (e_i(t))
 \end{equation}
for each $i=1,\cdots,n$, where $\mathcal{W}=(h_i^j)$ denotes the Weingarten map.  The resulting local frame field $\{e_i(t)\}_{i=1}^n$ remains orthonormal with respect to $g(x,t)$.  In fact, the evolution equation \eqref{s3:2-gij} for the metric implies that
\begin{align*}
  \frac {\partial}{\partial t}\left(g(e_i(t),e_j(t))\right)=&(\frac {\partial}{\partial t}g)(e_i(t),e_j(t))+g( \frac {\partial}{\partial t}e_i(t),e_j(t))+g(e_i(t), \frac {\partial}{\partial t}e_j(t))\\
  =&2\left(\frac{\lambda'(r)-u}{F(\tilde{\kappa})}-u\right) h(e_i(t),e_j(t))\\
  &-\left(\frac{\lambda'(r)-u}{F(\tilde{\kappa})}-u\right) \left(g(\mathcal{W}(e_i(t)),e_j(t)) +g(e_i(t), \mathcal{W} (e_j(t)))\right)\\
  =& 0.
\end{align*}
We now calculate the evolution of $\zeta$:
\begin{align*}
  \frac {\partial}{\partial t} \zeta =&  \frac {\partial}{\partial t} \left(S(e_n(t),e_n(t))\right) \\
  = &  \left(\frac{\partial}{\partial t} S\right)(e_n(t),e_n(t)) +2 S\left(\frac{\partial}{\partial t}e_n(t),e_n(t)\right) \\
  = &   \left(\frac {\partial}{\partial t} S\right) (e_n(t),e_n(t)) -2 \left(\frac{\lambda'(r)-u}{F(\tilde{\kappa})}-u\right)S(\mathcal{W}(e_n(t)), e_n(t))\\
  =&  \frac {\partial}{\partial t} S_{nn}-2 \left(\frac{\lambda'(r)-u}{F(\tilde{\kappa})}-u\right)S_{nn}h_{nn},
\end{align*}
where in the third equality we used \eqref{s7.3pf1} for the evolution of $e_n(t)$. Substituting the evolution equation \eqref{s3:2-hij} for $S_{ij}$ and setting $i=j=n$,  we have
	\begin{align}\label{s7.zeta}
	\frac{\partial}{\partial t} \zeta=&\frac{\lambda'-u}{F^2}\dot{F}^{kl}\nabla_k\nabla_l\zeta+\frac{\lambda'-u}{F^2}\ddot{F}^{kl,pq}\nabla_nS_{kl}\nabla_nS_{pq}+\frac 2F\nabla_nF\nabla_n(\frac{\lambda'-u}F)\nonumber\\
&+\frac {1+F}F\langle V,\nabla \zeta \rangle+\left(\lambda'+\frac{\lambda'-u}{F^2}(\dot{F}^{kl}(S^2)_{kl}-F)\right)\zeta\nonumber\\
&-\left(\frac{\lambda'-u}{F^2}(F+\dot{F}^{kl}g_{kl})+\frac{\lambda'}F\right)\zeta^2+(\lambda' -u) \left(1+\frac{\dot{F}^{kl}(S^2)_{kl}}{F^2}\right).
	\end{align}
Keep in mind that $\zeta=S_{nn}$ and $h_{nn}=\zeta+1$.  The second term on the right hand side of \eqref{s7.zeta} is non-positive by the concavity of $F$. The third term  of \eqref{s7.zeta} can be estimated as
\begin{align*}
 \frac 2F\nabla_nF\nabla_n(\frac{\lambda'-u}F)=& -\frac {2(\lambda'-u)}{F^3}|\nabla_nF|^2+\frac 2{F^2}\nabla_nF\nabla_n(\lambda'-u) \\
 =& -\frac {2(\lambda'-u)}{F^3}\left(\nabla_nF-\frac{F}{2(\lambda'-u)}\nabla_n(\lambda'-u)\right)^2 +\frac{|\nabla_n(\lambda'-u)|^2}{2(\lambda'-u)F}\\
 \leq   &~ \frac{\lambda^2-u^2}{2(\lambda'-u)F}\zeta^2.
\end{align*}
Using \eqref{s2:0-1} - \eqref{s2:0-2} and noting that $2\lambda'(\lambda'-u)>\lambda^2-u^2$, we have
\begin{align}\label{s7.zeta2}
	 &\frac \partial{\partial t} \zeta-\frac{\lambda'-u}{F^2}\dot{F}^{kl}\nabla_k\nabla_l\zeta-\frac {1+F}F\langle V,\nabla \zeta \rangle \nonumber\\
\leq &\left(\lambda'+\frac{\lambda'-u}{F^2}(\dot{F}^{kl}(S^2)_{kl}-F)\right)\zeta\nonumber\\
&-\left(\frac{\lambda'-u}{F^2}(F+\dot{F}^{kl}g_{kl})+\frac{\lambda'}F-\frac{\lambda^2-u^2}{2(\lambda'-u)F}\right)\zeta^2\nonumber\\
 &+(\lambda' -u) \left(1+\frac{\dot{F}^{kl}(S^2)_{kl}}{F^2}\right)\nonumber\\
 \leq &\left(\lambda'+\frac{\lambda'-u}{F^2}((n-m+1)F^2-F)\right)\zeta\nonumber\\
 &-\frac{\lambda'-u}{F^2}(F+1)\zeta^2+(n-m+2)(\lambda' -u).
	\end{align}
By the $C^0$ estimate and the bounds on $F$, the dominated term on the right hand side of of \eqref{s7.zeta2} is the term involving $\zeta^2$. We can apply the maximum principle to conclude that $\zeta$ is bounded from above.
\endproof

\section{Long time existence and convergence}\label{sec:conver}
In this section, we prove the long time existence and exponential convergence of the flow \eqref{s1:flow1}, and complete the proof of Theorem \ref{main-thm-II}.

\subsection{Long time existence}$\ $

\begin{prop}\label{pro-gradient-estimate}
Let $\g=\g(\t,t)$, $(\t,t)\in \mathbb S^n \times [0,T)$ be a solution of the initial value problem \eqref{s6:varphi}, then
	\begin{align}\label{C1-estimate}
	|D \g|^2(\t,t) \leq C, \quad \forall (\t,t)\in \mathbb{S}^n \times [0,T),
	\end{align}
	where $C>0$ is a constant depending only on the initial hypersurface $M_0$.
\end{prop}
\proof
Note that $v=\sqrt{1+|D\varphi|^2}$ by the definition \eqref{s2:v-def}, the estimate \eqref{C1-estimate} is equivalent to an upper bound on $v$. Since the support function $u=\lambda/v$, by the $C^0$ estimate, it suffices to prove a positive lower bound on $u$.  For each positive time $t$, at the spatial minimum point of the support function $u$, we have
	$$
	h_i^k\nabla_k \l'=\nabla_i u=0.
	$$
	Since the hypersurface $M_t$ is strictly h-convex, we have
\begin{equation*}
  0=\nabla_i\lambda'(r)=\langle \lambda(r)\partial_r,e_i\rangle,\quad i=1,\cdots,n ,
\end{equation*}
at the spatial minimum point $P$ of $u$, where the second equality is due to \eqref{2.1}. Equivalently, the unit normal vector $\nu$ is proportional to $\partial_r$ at $P$. Thus, $u|_{P}=\l|_{P}$ and hence
	\begin{align*}
	\frac{\l(r)}{v}=u \geq \min_{M_t} u \geq \min_{M_t} \l(r)>0.
	\end{align*}
	This implies that
	$$
	v\leq \frac{\max_{M_t}\l}{\min_{M_t}\l}\leq \frac{\max_{M_0}\l}{\min_{M_0}\l} \leq C
	$$
by the $C^0$ estimate. From the definition of $v$, we get the gradient estimate \eqref{C1-estimate}.
\endproof

By Propositions \ref{pro-C0-estimate} and \ref{pro-gradient-estimate}, we have the $C^0$ and $C^1$-estimates of the solution $\g$. On the other hand, by Theorem \ref{thm-h-convexity-2} and Proposition \ref{s8:prop-upper-bound}, we have $0< \~\k_i \leq C$, $1\leq i\leq n$. Together with $C^0$ and $C^1$-estimates, the curvature estimates give the uniform $C^2$-estimate of the solution $M_t$. Moreover, the uniform positive lower bound on $F(\tilde{\kappa})$ and the upper bound on $\tilde{\kappa}_i$ imply that $\tilde{\kappa}_i$ lies in a compact subset of the Garding cone $\Gamma_m^+$. Therefore, the flow equation \eqref{s1:flow1} is uniformly parabolic with concave operator in the second spatial derivatives. Now we can apply Krylov's \cite{Kryl82} theory to derive uniform $C^{2,\a}$-estimate of the solution $M_t$, and then apply standard parabolic Schauder estimate to derive the higher regularity estimate. Consequently, we obtain the long time existence of the flow \eqref{s1:flow1}.

\subsection{Convergence to a geodesic sphere}\label{sec:umbi}$\ $

Along the flow \eqref{s1:flow1}, we have
	\begin{align*}
	\frac{d}{dt}\widetilde{W}_{m+1}(\Omega_t)=&~\int_{M_t}\left( (\l'-u)\frac{E_{m+1}(\~\k)E_{m-1}(\~\k)}{E_{m}(\~\k)}-uE_{m+1}(\~\k)\right) d\mu_t  \\
	= &~\int_{M_t}(\l'-u)\left( \frac{E_{m+1}(\~\k)E_{m-1}(\~\k)}{E_{m}(\~\k)}-E_{m}(\~\k)\right) d\mu_t ~\leq~0,
	\end{align*}
where the last inequality follows from the Newton-MacLaurin inequality \eqref{s2:Newt-2}. Furthermore, by the $C^0$ estimate, each $\Omega_t$ contains a geodesic ball of radius $r_{\min}(0)>0$, where $r_{\min}(0)$ is the minimum of the radial function of the initial hypersurface $M_0=\partial\Omega_0$. By the monotonicity of $\widetilde{W}_{m+1}$ under the inclusion of h-convex domains (\cite[Corollary 5.2]{ACW2018}), we have
\begin{equation*}
  \widetilde{W}_{m+1}(\Omega_t)\geq~\widetilde{W}_{m+1}(B_{r_{\min}(0)})>0.
\end{equation*}
Then we have
\begin{equation*}
 0\leq~ -\int_0^\infty\frac{d}{dt}\widetilde{W}_{m+1}(\Omega_{t}) dt ~<\infty.
\end{equation*}
Since $\l'-u\geq C>0$, by the regularity estimate of the solution $M_t$, for any $\varepsilon>0$, there exists a $T_0>0$ such that for all $t>T_0$, there holds
\begin{equation*}
 0\leq E_{m}(\~\k)^2-E_{m+1}(\~\k)E_{m-1}(\~\k)\leq \varepsilon.
\end{equation*}
Equivalently,
\begin{equation}\label{s8:2-umb}
  \max_{i<j}|\kappa_i-\kappa_j|\leq~\varepsilon
\end{equation}
holds on $M_t$ for all $t>T_0$.

On the other hand, the regularity estimate of the solution implies that there exists a sequence of times $\{t_i\}$ such that $\Omega_{t_i}$ converges to $\Omega_{\infty}$ smoothly as $t_i\to\infty$. By the above analysis, $M_{\infty}=\partial\Omega_{\infty}$ must be a geodesic sphere. Since the velocity of the flow \eqref{s1:flow1} has no global terms, by the comparison principle, all the evolving hypersurfaces $M_t$ must converge to the unique geodesic sphere $M_{\infty}$.

The limit geodesic sphere must be centered at the origin. If not, by the proof of Proposition \ref{pro-C0-estimate},
\begin{equation*}
  \frac d{dt}r_{\max}<0,\quad \mathrm{and}\quad \frac d{dt}r_{\min}>0
\end{equation*}
on $M_{\infty}$. The flow can immediately move the geodesic sphere away. Hence we conclude that the solution of \eqref{s1:flow1} converges smoothly to a unique geodesic sphere $S_{r_\infty}(0)$ centered at the origin as time $t\to\infty$. The radius $r_\infty$ can be determined by the preservation of $\widetilde{W}_m(\Omega_t)$, i.e., the radius $r_\infty$ is the one such that
\begin{equation*}
  \widetilde{W}_m(B_{r_\infty})=\widetilde{W}_m(\Omega_0).
\end{equation*}

\subsection{Exponential convergence}$\ $

Finally, we show that the convergence is in exponentially rate. The convergence to a geodesic sphere centered at the origin implies that $|D\varphi|$ decays to zero as $t\to\infty$. We next prove that $|D\varphi|$ decays to zero exponentially.

Define $\omega=\frac 12|D\varphi|^2$. We first calculate the evolution equation of $\omega$ using the evolution equation \eqref{s6:varphi} of $\varphi$.
\begin{lem}
	The function $\omega=\frac 12|D\varphi|^2$ satisfies the following equation along the flow \eqref{s1:flow1}
	\begin{align}\label{s4:2-omeg}
	\frac{\partial}{\partial t}\omega =&\frac{2\lambda v\omega}{ F}-\frac{2\lambda'\omega}{F}+\frac{2\lambda \omega}{vF^2} \left(\frac{\lambda'v}{\lambda}-1\right)^2\dot{F}^{ij}\delta_i^j\nonumber\\
	&+\frac{\lambda'}{\lambda vF}\varphi^i \omega_i+\frac{3}{v^2F^2}\left(\frac{\lambda'v}{\lambda}-1\right)\varphi^k \omega_k(F+\dot{F}^{ij}\delta_i^j)\nonumber\\
	&-\frac 2{\lambda v^3F^2}\left(\frac{\lambda'v}{\lambda}-1\right)\varphi^k \omega_k\dot{F}^{ij}\left(-\varphi_i^j+\lambda'\delta_i^j\right)\nonumber\\
	&-\frac 1{\lambda vF^2}\left(\frac{\lambda'v}{\lambda}-1\right)\dot{F}^{ij}\left(\varphi_{i}^{k}\varphi_{k}^j-\varphi_i\varphi^j+|D\varphi|^2\delta_i^j\right)\nonumber\\
	& -\frac 1{\lambda v^3F^2}\left(\frac{\lambda'v}{\lambda}-1\right)\dot{F}^{ij}\left(\omega_i\omega^j +(\varphi_k\varphi^j-v^2\delta_k^j) \omega_{i}^{k}\right).
	\end{align}
\end{lem}
\proof
Since the support function $u=\lambda/v$, where $v=\sqrt{1+|D \g|^2}$, the equation \eqref{s6:varphi} is equivalent to
\begin{equation}\label{s4:2-1}
\frac{\partial}{\partial t} \g =\left(\frac{\lambda'v}{\lambda}-1\right)\frac 1F-1.
\end{equation}
We have
\begin{align*}
\frac{\partial}{\partial t}\omega= & \varphi^i D_i(\frac{\partial}{\partial t}\varphi) \\
= & \frac 1F \varphi^i \left(\frac 1{\lambda}(vD_i\lambda'+\lambda' v_i)-\frac{\lambda'v}{\lambda^2}D_i\lambda\right)-\frac 1{F^2}\left(\frac{\lambda'v}{\lambda}-1\right)\varphi^i D_iF
\end{align*}
Note that
\begin{align*}
D_k\lambda= & \lambda' r_k=\lambda\lambda'\varphi_k \\
D_k\lambda'= & \lambda''r_k=\lambda^2\varphi_k\\
v_k=&\frac 1v\omega_k
\end{align*}
Then
\begin{align}\label{s4:2-2}
\frac{\partial}{\partial t}\omega  = & \frac 1F \varphi^i \left(\frac 1{\lambda}(\lambda^2 v\varphi_i+\frac{\lambda'}v \omega_i)-\frac{(\lambda')^2v}{\lambda}\varphi_i\right)-\frac 1{F^2}\left(\frac{\lambda'v}{\lambda}-1\right)\varphi^i D_iF\nonumber\\
=&-\frac{2v\omega}{\lambda F}+\frac{\lambda'}{\lambda vF}\varphi^i \omega_i-\frac 1{F^2}\left(\frac{\lambda'v}{\lambda}-1\right)\varphi^i D_iF.
\end{align}
We next calculate the term $\varphi^i D_iF$: By \eqref{s2:h}, the Weingarten matrix $h_i^j$ can be expressed as following
\begin{align}\label{s8:h}
h_i^j=&-\frac{1}{\l v}(e^{jk}-\frac{\g^j\g^k}{v^2})\varphi_{ik}+\frac{\lambda'}{\lambda v}\delta_i^j \nonumber\\
=& \frac{1}{\l v^3}\biggl(v^2\left(-\varphi_i^j+\lambda'\delta_i^j\right)+\varphi^j \omega_i\biggr),
\end{align}
where $\varphi_i^j=e^{jk}\varphi_{ki}$. We have
\begin{align*}
\varphi^k D_kF =& \varphi^k \dot{F}^{ij}D_kh_i^j \\
= & \varphi^k \dot{F}^{ij}D_k( \frac{1}{\l v^3}) \lambda v^3h_i^j\\
&+\frac{1}{\l v^3}\varphi^k \dot{F}^{ij}D_k\biggl(v^2\left(-\varphi_i^j+\lambda'\delta_i^j\right)+\varphi^j \omega_i\biggr) \\
=& -\frac{F+\dot{F}^{ij}\delta_i^j}{\l v^3} \varphi^k (v^3 D_k\lambda +3\lambda v^2 v_k)+\frac{2}{\l v^2}\varphi^k v_k\dot{F}^{ij}\left(-\varphi_i^j+\lambda'\delta_i^j\right)\\
&+\frac{1}{\l v^3}\varphi^k \dot{F}^{ij}\biggl(v^2\left(-\varphi_i{}^j{}_k+D_k\lambda'\delta_i^j\right)+\varphi^j_k \omega_i+\varphi^j \omega_{ki}\biggr)\\
=& -2\lambda'(F+\dot{F}^{ij}\delta_i^j)\omega-\frac{3(F+\dot{F}^{ij}\delta_i^j)}{v^2}\varphi^k \omega_k\\
&+\frac{2}{\l v^3}\varphi^k \omega_k\dot{F}^{ij}\left(-\varphi_i^j+\lambda'\delta_i^j\right)\\
&+\frac{1}{\l v^3}\varphi^k \dot{F}^{ij}\biggl(v^2\left(-\varphi_i{}^j{}_k+\lambda^2\varphi_k\delta_i^j\right)+\varphi^j_k\omega_i+\varphi^j \omega_{ki}\biggr).
\end{align*}
Note that
\begin{align*}
\omega_i=&\varphi^k\varphi_{ki}\\
\omega_{ij}=&\varphi^k_j \varphi_{ki}+\varphi^k \varphi_{kij}\\
=& \varphi^{k}_j\varphi_{ki}+\varphi^k\left(\varphi_{ijk}-\varphi_j\delta_k^i+\varphi_k\delta_i^j\right)\\
=& \varphi^{k}_j\varphi_{ki}+\varphi^k \varphi_{ijk}-\varphi_i\varphi_j+|D\varphi|^2\delta_i^j,
\end{align*}
where we used the Ricci identity to commute the covariant derivatives on the sphere $\mathbb{S}^n$. We have
\begin{align}\label{s4:2-3}
\varphi^k D_kF
=& -2\lambda'F\omega-\frac{2\lambda\omega}v \left(\frac{\lambda'v}{\lambda}-1\right)\dot{F}^{ij}\delta_i^j-\frac{3(F+\dot{F}^{ij}\delta_i^j)}{v^2}\varphi^k \omega_k\nonumber\\
&+\frac{2}{\l v^3}\varphi^k \omega_k\dot{F}^{ij}\left(-\varphi_i^j+\lambda'\delta_i^j\right)\nonumber\\
&+\frac{1}{\l v}\dot{F}^{ij}\left(-\omega_i^j+\varphi^{k}_{i}\varphi_{k}^{j}-\varphi_i\varphi^j+|D\varphi|^2\delta_i^j\right)\nonumber\\
&+\frac{1}{\l v^3}\dot{F}^{ij}\left(\omega^j \omega_i+\varphi^k \varphi_j \omega_{ki}\right).
\end{align}
Substituting \eqref{s4:2-3} into \eqref{s4:2-2} and rearranging terms, we obtain the equation \eqref{s4:2-omeg}.
\endproof

\begin{prop}\label{pro-exponential-convergence}
	Let $\g=\g(\t,t)$, $(\t,t)\in \mathbb S^n \times [0,\infty)$ be a solution of the initial value problem \eqref{s6:varphi}, and $r_\infty$ be the radius of the limit geodesic sphere as given in Section \ref{sec:umbi}. For any positive constant
 \begin{equation*}
   \alpha<\frac{2(n-1)}{n(\lambda'(r_\infty)-\lambda(r_\infty))}
 \end{equation*}
which depends only on the initial hypersurface $M_0$,  we have
	\begin{align}\label{exponential-convergence}
	|D \g|^2(\t,t) \leq C e^{-\a t}, \quad \forall (\t,t)\in \mathbb{S}^n \times [0,\infty),
	\end{align}
where $C$ is a positive constant depending on $M_0$ and $\alpha$.
\end{prop}
\begin{proof}
Since we focus on the limit behavior of the flow, we only need to get the estimate of $\omega=\frac 12|D \varphi|^2$ for time $t>T_0$, where $T_0$ is the time in Section \ref{sec:umbi} such that the almost total umbilicity \eqref{s8:2-umb} holds on $M_t$. Without loss of generality, we also assume that $\omega\leq \varepsilon$ and $\dot{F}^{ii}\in (1/n-\varepsilon,1/n+\varepsilon)$ hold on $M_t$ for time $t>T_0$ (by possibly increasing $T_0$).

For any time $t_0>T_0$, let $\theta_{t_0}\in \mathbb{S}^n$ be the maximum point of $\omega$ at $t=t_0$. Then at $(\theta_{t_0},t_0)$, we have
\begin{align*}
D  \omega = 0, \quad  (D^2 \omega) \leq 0.
\end{align*}
Substituting these into \eqref{s4:2-omeg}, we get
	\begin{align}\label{s8:3-exp1}
	\frac{\partial}{\partial t}\omega \leq &\frac{2\lambda v\omega}{ F}-\frac{2\lambda'\omega}{F}+\frac{2\lambda \omega}{vF^2} \left(\frac{\lambda'v}{\lambda}-1\right)^2\dot{F}^{ij}\delta_i^j\nonumber\\
	&-\frac 1{\lambda vF^2}\left(\frac{\lambda'v}{\lambda}-1\right)\dot{F}^{ij}\left(|D\varphi|^2\delta_i^j-\varphi_i\varphi^j\right).
	\end{align}
at $(\theta_{t_0},t_0)$, where we also used the fact that both $\dot{F}^{ij}$ and $(v^2\d_k^j-\varphi_k\varphi^j)$ are positive definite. At this point, we choose orthonormal frame such that
	    $$
	    \varphi_1=|D \varphi|,\quad \varphi_j=0, \quad 2\leq j\leq n.
	    $$
and $(\varphi_{ij})$ is diagonal with $\varphi_{1j}=0$ for all $j=1,\cdots,n$. Then by \eqref{s2:h}, the Weingarten matrix
\begin{equation}\label{s8:3-h}
  h_i^j=\frac 1{\lambda v}\left(-\varphi_{ij}+\lambda'\delta_i^j\right)
\end{equation}
is diagonalized at $(\theta_{t_0},t_0)$. By \eqref{s8:2-umb} and $\varphi_{11}=0$, we have $|\varphi_{ii}|=|\varphi_{ii}-\varphi_{11}|\leq \varepsilon \lambda v$.
By \eqref{s8:3-h} and the homogeneity of $F$,
\begin{align*}
  \lambda vF= &\sum_{i=1}^n \dot{F}^{ii}(-\varphi_{ii}+\lambda'-\lambda v) \\
   =&( \lambda'-\lambda v)\sum_{i=1}^n \dot{F}^{ii}-\sum_{i=1}^n \dot{F}^{ii}\varphi_{ii}
\end{align*}
The inequality \eqref{s8:3-exp1} becomes
 \begin{align*}
	\frac{\partial}{\partial t}\omega \leq &\frac{2\lambda v\omega}{ F}-\frac{2\lambda'\omega}{F}+\frac{2\lambda \omega}{vF^2} \left(\frac{\lambda'v}{\lambda}-1\right)^2\sum_{i=1}^n\dot{F}^{ii}\nonumber\\
	&-\frac {2\omega}{\lambda vF^2}\left(\frac{\lambda'v}{\lambda}-1\right)\sum_{i=2}^n\dot{F}^{ii}\nonumber\\
\leq &\frac{2\omega}{F^2}\biggl(\frac 1{\lambda v}\left((\lambda'v-\lambda)^2-(\lambda'-\lambda v)^2\right)\sum_{i=1}^n\dot{F}^{ii}\nonumber\\
	&-\frac {1}{\lambda^2 v}(\lambda'v-\lambda)\sum_{i=2}^n\dot{F}^{ii}-\frac{\lambda v-\lambda'}{\lambda v}\sum_{i=1}^n \dot{F}^{ii}\varphi_{ii}\biggr)\nonumber\\
\leq &-\left(\frac{2(n-1)}{n(\lambda'-\lambda)}-\tilde{\varepsilon}\right)\omega,
	\end{align*}
where $\tilde{\varepsilon}\leq C\varepsilon$ is a small constant. Let
\begin{equation*}
  \alpha_0=\frac{2(n-1)}{n(\lambda'(r_\infty)-\lambda(r_\infty))}.
\end{equation*}
We conclude that for any $\alpha<\alpha_0$, there exists a time $T_0$ depending on $M_0$ and $\alpha_0-\alpha$ such that for all $t>T_0$, there holds
\begin{equation*}
\frac{\partial}{\partial t}\omega \leq -\alpha \omega,\qquad t\in [T_0,\infty)
\end{equation*}
at the spatial maximum point of $\omega$. Integrating the above inequality, we obtain
\begin{equation*}
  \omega(t)\leq \omega(T_0)e^{-\alpha(t-T_0)}
\end{equation*}
for all time $t\geq T_0$. This gives the exponential decay \eqref{exponential-convergence} of $|D\varphi|^2$.
\end{proof}

\section{New geometric inequalities: II}\label{sec:9}
In this section, we apply the convergence result of the flow \eqref{s1:flow1} to prove the new geometric inequalities in Theorem \ref{thm-shift-weighted-AF-inequality}, Corollary \ref{corollary-shift-weighted-AF-inequality} and Corollary \ref{corollary-weighted-Lk-inequality}.

\subsection{Proof of Theorem \ref{thm-shift-weighted-AF-inequality}} We first calculate the evolution of $\int_{M_t}(\l'-u) E_k(\~\k)d\mu_t$ along the new constrained curvature flow \eqref{s1:flow1}. Along the general flow
\begin{equation*}
  \frac{\partial}{\partial t} X(x,t)= \eta(x,t)\nu(x,t)
\end{equation*}
in hyperbolic space, the support function $u$ evolves by
\begin{align}\label{s9:u-evl}
   \frac{\partial}{\partial t}u = & \langle  \frac{\partial}{\partial t} (\lambda(r)\partial_r),\nu\rangle+\langle \lambda(r)\partial_r, \frac{\partial}{\partial t} \nu\rangle\nonumber \\
  = &\lambda'(r)\eta-\langle\lambda(r)\partial_r,\nabla\eta\rangle \nonumber\\
  = &\lambda'(r)\eta-\langle\nabla\lambda'(r),\nabla\eta\rangle,
\end{align}
where we used the fact that $\lambda(r)\partial_r$ is a conformal Killing vector field. Combining \eqref{s5:dmu} -- \eqref{s5:2Ek} and \eqref{s9:u-evl}, we have
\begin{align*}
\frac{d}{dt}\int_{M_t}(\l'-u) E_k(\~\k) d\mu_t=&\int_{M_t}(\l'-u) E_k(\~\k) H\eta d\mu_t+\int_{M_t} u E_k(\~\k)\eta d\mu_t\\
&-\int_{M_t} (\l' \eta-\langle \nabla \l',\nabla \eta\rangle)E_k(\~\k)d\mu_t\\
& +\int_{M_t}(\l'-u) \dot{E}_k^{ij} \left(-\nabla_i\nabla^j \eta-\eta((h^2)_i^j-\d_i^j)\right)d\mu_t,
\end{align*}
where $\dot{E}_k^{ij} $ is evaluated at the shifted Weingarten matrix $(S_i^j)=(S_{ik}g^{kj})$. Since $S_{ij}=h_{ij}-g_{ij}$ is Codazzi, the operator $\dot{E}_k^{ij} $ is divergence-free. By integration in part, we have
\begin{align*}
\frac{d}{dt}\int_{M_t}(\l'-u)E_k(\~\k) d\mu_t=&\int_{M_t} (\l'-u)\left(-E_k(\~\k)-\dot{E}_k^{ij}((h^2)_i^j-\d_i^j)+E_k(\~\k)H\right) \eta  d\mu_t\\
                                              &+\int_{M_t} \langle \nabla \l',\nabla \eta\rangle E_k(\~\k) d\mu_t-\int_{M_t}\nabla_i \nabla^j(\l'-u) \dot{E}_k^{ij} \eta  d\mu_t\\
=&\int_{M_t} (\l'-u)\left(-E_k(\~\k)-\dot{E}_k^{ij}((h^2)_i^j-\d_i^j)+E_k(\~\k)H\right) \eta  d\mu_t\\
&+\int_{M_t} \langle \nabla \l',\nabla \eta\rangle E_k(\~\k) d\mu_t-\int_{M_t}\dot{E}_k^{ij}(\l'\d_i^j-u h_i^j)\eta  d\mu_t\\
&+\int_{M_t} \dot{E}_k^{ij} \left(\langle \nabla \l',\nabla h_i^j\rangle+\l' h_i^j-u(h^2)_i^j\right)  \eta  d\mu_t\\
=&\int_{M_t}\l'\biggl((n-k)E_{k+1}(\~\k)+(n-k-1)E_k(\~\k)\biggr)\eta  d\mu_t\\
 & +\int_{M_t} u \left( (k+1) E_k(\~\k)-E_k(\~\k)H \right) \eta+\int_{M_t} \langle \nabla \l',\nabla (\eta E_k(\~\k))\rangle  d\mu_t,
\end{align*}
where we used \eqref{newton-formula-1} -- \eqref{newton-formula-3}, \eqref{2.1} -- \eqref{2.2} and $H=nE_1(\~\k)+n$. By divergence theorem and \eqref{2.1}, we have
\begin{align*}
\int_{M_t} \langle \nabla \l',\nabla (\eta E_k(\~\k))\rangle  d\mu_t=-\int_{M_t} (n\l'-Hu)E_k(\~\k)\eta  d\mu_t.
\end{align*}
This gives us
\begin{align*}
\frac{d}{dt}\int_{M_t}(\l'-u)E_k(\~\k) d\mu_t=&\int_{M_t} \biggl(u(k+1) E_k(\~\k) + \l'\((n-k)E_{k+1}(\~\k)-(k+1)E_k(\~\k)\)\biggr)\eta  d\mu_t \\
                                             =&\int_{M_t} \biggl( (n-k)\l'E_{k+1}(\~\k)- (k+1)(\l'-u)E_{k}(\~\k)\biggr)\eta  d\mu_t.
\end{align*}

If $k=1,\cdots,n$ and $\~\kappa\in \G_k^{+}$, we choose
\begin{align}\label{s9:flow2}
\eta= \frac{(\l'-u)E_{k-1}(\~\k)}{E_k(\~\k)}-u;
\end{align}
while if $k=1,\cdots,n-1$ and $\~\kappa\in \G_{k+1}^{+}$, we choose
\begin{align}\label{s9:flow3}
\eta= \frac{(\l'-u)E_{k}(\~\k)}{E_{k+1}(\~\k)}-u.
\end{align}
Using the Newton-MacLaurin inequality \eqref{s2:Newt-2}, for \eqref{s9:flow2} and \eqref{s9:flow3} we always have
\begin{align*}
\frac{d}{dt}\int_{M_t} (\l'-u)E_k(\~\k) d\mu_t
\leq &(n-k)\int_{M_t} \l'\biggl( (\l'-u) E_{k}(\~\k)-uE_{k+1}(\~\k) \biggr) d\mu_t\\
     &-(k+1)\int_{M_t} (\l'-u)\biggl( (\l'-u) E_{k-1}(\~\k)-uE_{k}(\~\k) \biggr) d\mu_t.
\end{align*}
As in \eqref{s5:2-1}, for each $k=1,\cdots,n$ we have
\begin{align*}
\dot{E}_k^{ij} \nabla_i\nabla_j \l'=&\dot{E}_k^{ij}((\l'-u)\d_i^j-uS_i^j) \\
                                   =&k ((\l'-u) E_{k-1}(\~\k)-uE_k(\~\k)).
\end{align*}
Since $\dot{E}_k^{ij} $ and $\dot{E}_{k+1}^{ij} $ are divergence-free, we obtain
\begin{align}\label{s9:2}
\frac{d}{dt}\int_{M_t} (\l'-u)E_k(\~\k) d\mu_t \leq &~\frac{n-k}{k+1}\int_{M_t} \l' \dot{E}_{k+1}^{ij} \nabla_i \nabla_j \l'  d\mu_t-\frac{k+1}{k}\int_{M_t} (\l'-u)\dot{E}_{k}^{ij} \nabla_i \nabla_j \l'  d\mu_t\nonumber\\
=&-\frac{n-k}{k+1}\int_{M_t}\dot{E}_{k+1}^{ij} \nabla_i\lambda' \nabla_j \l' d\mu_t+\frac{k+1}{k}\int_{M_t} \dot{E}_{k}^{ij} \nabla_i (\l'-u)\nabla_j \l'  d\mu_t\nonumber\\
= &-\frac{n-k}{k+1}\int_{M_t} \dot{E}_{k+1}^{ij} \nabla_i \l'\nabla_j \l' d\mu_t-\frac{k+1}{k}\int_{M_t}\dot{E}_k^{ij}S_i^l \nabla_l \l' \nabla_j \l' d\mu_t,
\end{align}
where we used $\nabla_i (\l'-u)=-(h_i^l-\d_i^l)\nabla_l \l'=-S_i^l\nabla_l \l'$. By the strict h-convexity of $M_t$, all the matrices $\dot{E}_{k+1}^{ij}$, $\dot{E}_k^{ij}$ and $S_j^l$ are positive definite. This gives
\begin{align*}
\frac{d}{dt}\int_{M_t} (\l'-u)E_k(\~\k) d\mu_t \leq 0.
\end{align*}
Since the flow with speed given by \eqref{s9:flow2} preserves the new $k$th quermassintegral $\widetilde{W}_k(\Omega_t)$ and converges to a geodesic sphere of radius $r_{\infty}$, we obtain
\begin{equation}\label{s9:3}
  \int_{M} (\lambda'-u)E_k(\~\k) d\mu\geq \~h_k \circ \~f_k^{-1}(\widetilde{W}_k(\Omega)), \quad k=1,\cdots,n,
\end{equation}
where $\~h_k(r)=\omega_{n} (\cosh r-\sinh r)^{k+1}\sinh^{n-k}r$. Similarly, the flow with speed given by \eqref{s9:flow3} preserves $\widetilde{W}_{k+1}(\Omega_t)$ and converges to a geodesic sphere of radius $r_{\infty}$, we obtain
\begin{equation}\label{s9:4}
\int_{M} (\lambda'-u)E_k(\~\k) d\mu\geq \~h_k \circ \~f_{k+1}^{-1}(\widetilde{W}_{k+1}(\Omega)), \quad k=1,\cdots,n-1.
\end{equation}
If equality holds in \eqref{s9:3} or \eqref{s9:4}, the inequality \eqref{s9:2} implies that $\nabla \l'(r)=0$ on $M_t$ for all $t$. This implies that the initial hypersurface is a geodesic sphere centered at the origin.

\subsection{Proof of Corollary \ref{corollary-shift-weighted-AF-inequality}}

In general, the function $\tilde{h}_k(r)$ is not increasing in $r$. For $k \leq \frac{n-1}{2}$, we have
	$$
	\~h'_k(r)=-\omega_{n} e^{-(1+k)r} \sinh^{n-k-1}r \((k+1)\sinh r-(n-k)\cosh r\)>0.
	$$
Hence $\~h_k$ is strictly increasing in $r$ for $1 \leq k \leq \frac{n-1}{2}$. This implies the Corollary \ref{corollary-shift-weighted-AF-inequality}. In fact, let $1 \leq k \leq \frac{n-1}{2}$. If $\~\k\in \G_k^{+}$, then by \eqref{s1:ACW} and \eqref{s9:3} we have
	\begin{align*}
	\int_{M} (\l'-u)E_k(\~\k) d\mu ~\geq~ \~h_k \circ \tilde{f}_{k}^{-1}(\widetilde{W}_{k}(\Omega)) \geq \~h_k \circ \~f_\ell^{-1}(\widetilde{W}_\ell(\Omega)), \quad 0\leq \ell \leq k.
	\end{align*}
	Similarly, if $\~\kappa \in \G_{k+1}^{+}$, then by \eqref{s1:ACW} and \eqref{s9:4} we have
	\begin{align*}
	\int_{M} (\l'-u)E_k(\~\k) d\mu ~\geq~ \~h_k \circ \tilde{f}_{k+1}^{-1}(\widetilde{W}_{k+1}(\Omega)) \geq \~h_k \circ \~f_\ell^{-1}(\widetilde{W}_\ell(\Omega)), \quad 0\leq \ell \leq k+1.
	\end{align*}
Both inequalities are sharp with equality holds if and only if $M$ is a geodesic sphere centered at the origin.

\subsection{Proof of Corollary \ref{corollary-weighted-Lk-inequality}}

 We first show that the  Gauss-Bonnet curvature $L_k$ can be expressed as a linear combination of the shifted $m$th mean curvatures with $m$ ranging from $k$ to $2k$.
\begin{lem}\label{s9:lem-Lk}
	For a hypersurface $(M,g)$ in $\mathbb H^{n+1}$, its Gauss-Bonnet curvature $L_k$ of the induced metric $g$ can be expressed as follows
	\begin{align}\label{s9:Lk-shifted-curvature-integral}
	L_k = \binom{n}{2k}(2k)!\sum_{j=0}^{k}2^j\binom kj E_{2k-j}(\~\k).
	\end{align}
\end{lem}
\begin{proof}
	Firstly, we recall the Gauss formula for hypersurfaces in hyperbolic space:
	\begin{align}\label{s9.Gauss}
	R_{ij}{}^{sl}=&(h_i^s h_j^l -h_i^l h_j^s)- (\d_i^s \d_j^l -\d_i^l \d_j^s)\nonumber\\
	=&(\~h_i^s+\d_i^s) (\~h_j^l+\d_j^l)-(\~h_i^l+\d_i^l)(\~h_j^s+\d_j^s)-(\d_i^s \d_j^l -\d_i^l \d_j^s) \nonumber\\
	=&(\~h_i^s\~h_j^l+\d_i^s \~h_j^l +\~h_i^s \d_j^l)-(\~h_i^l\~h_j^s +\d_i^l \~h_j^s +\~h_i^l \d_j^s),
	\end{align}
	where $\~h_i^j=h_i^j-\d_i^j$. Substituting \eqref{s9.Gauss} into the definition \eqref{Gauss-Bonnet-curvature} of $L_k$, we have
	\begin{align*}
	L_k
	=& \d_{j_1j_2\cdots j_{2k-1} j_{2k}}^{i_1i_2\cdots i_{2k-1}i_{2k}}(\~h_{i_1}^{j_1}\~h_{i_2}^{j_2}+2\~h_{i_1}^{j_1}\d_{i_2}^{j_2})\cdots (\~h_{i_{2k-1}}^{j_{2k-1}}\~h_{i_{2k}}^{j_{2k}}+2\~h_{i_{2k-1}}^{j_{2k-1}}\d_{i_{2k}}^{j_{2k}}) \\
	=& \sum_{m=0}^{k} \binom{k}{m} 2^m \d_{j_1 j_2 \cdots j_{2k}}^{i_1i_2\cdots i_{2k}} \d_{i_2}^{j_2}\d_{i_4}^{j_4}\cdots\d_{i_{2m}}^{j_{2m}} \~h_{i_1}^{j_1}\~h_{i_3}^{j_3}\cdots \~h_{i_{2m-1}}^{j_{2m-1}}\~h_{i_{2m+1}}^{j_{2m+1}}\~h_{i_{2m+2}}^{j_{2m+2}}\cdots \~h_{i_{2k-1}}^{j_{2k-1}}\~h_{i_{2k}}^{j_{2k}} \\
	=& \sum_{m=0}^{k} \binom{k}{m} 2^m (n+1-2k)\cdots (n+1-2k+m-1) (2k-m)! \binom{n}{2k-m}E_{2k-m}(\~\k)\\
	=& \binom{n}{2k}(2k)!\sum_{m=0}^{k}\binom{k}{m} 2^m E_{2k-m}(\~\k).
	\end{align*}
	Here in the first equality we used the symmetry of the generalized Kronecker delta, and in the third equality we used the basic property of the generalized Kronecker delta
	\begin{align*}
	\d_{j_1 j_2 \cdots j_p}^{i_1i_2\cdots i_{p}} \d_{i_1}^{j_1} =(n+1-p) \d_{j_2\cdots j_p}^{i_2\cdots i_p}.
	\end{align*}
\end{proof}	

\begin{proof}[Proof of Corollary \ref{corollary-weighted-Lk-inequality}]
	Using the expression \eqref{s9:Lk-shifted-curvature-integral} and shifted Minkowski formula \eqref{s2:shift-MF}, we have
	\begin{align}
	\int_M u L_k d\mu = & \binom{n}{2k} (2k)!  \sum_{j=0}^{k}2^j\binom kj \int_M uE_{2k-j}(\~\k) d\mu \nonumber\\
	                  = & \binom{n}{2k} (2k)! \sum_{j=0}^{k}2^j\binom kj\int_M (\l'-u) E_{2k-j-1}(\~\k) d\mu \nonumber\\
	                  \geq & \binom{n}{2k} (2k)!\sum_{j=0}^{k}2^j\binom kj \~h_{2k-j-1} \circ \~f_\ell^{-1}(\widetilde{W}_\ell(\Omega)) \nonumber\\
	                  = & \~g_k \circ \~f_\ell^{-1}(\widetilde{W}_\ell(\Omega)), \label{s9:Lk-shifted-quermassintegral}
	\end{align}
	where we used Corollary \ref{corollary-shift-weighted-AF-inequality} in the third inequality, provided that $2\leq k\leq \frac{n+1}{4}$, $\~\k\in\G_{2k-1}^{+}$ and $0\leq \ell<k$. The last equality follows from
	\begin{align*}
	\~g_k(r)=&\binom{n}{2k}(2k)!\sum_{j=0}^{k}2^j\binom kj \~h_{2k-j-1}(r)\\
	=&\binom{n}{2k}(2k)! \omega_n  (\cosh r-\sinh r)^k \sinh^{n+1-2k}r \sum_{j=0}^{k}\binom kj  (\cosh r-\sinh r)^{k-j}(2\sinh r)^j \\
	=&\binom{n}{2k}(2k)! \omega_n  \sinh^{n+1-2k}r.
	\end{align*}
	By the isoperimetric inequality $W_1(\Omega)\geq f_1\circ f_0^{-1}(W_0(\Omega))$ in hyperbolic space, we have
	\begin{align}
	\widetilde{W}_1(\Omega) = & W_1(\Omega)-W_0(\Omega) \nonumber\\
	                        \geq & W_1(\Omega)-f_0 \circ f_1^{-1}(W_1(\Omega)) \nonumber\\
	                        = & (f_1- f_0) \circ f_1^{-1}(W_1(\Omega)) \nonumber\\
	                        = & \~f_1 \circ f_1^{-1}(W_1(\Omega)).  \label{s9:isoperi-ineq}
	\end{align}
    Since $\~g_k$ is strictly increasing in $r$, it follows from \eqref{s9:Lk-shifted-quermassintegral} with $\ell=1$ and \eqref{s9:isoperi-ineq} that
	\begin{align}\label{s9:Lk-area}
	\int_M u L_k d\mu \geq \~g_k \circ \~f_1^{-1}(\widetilde{W}_1(\Omega)) \geq \~g_k \circ f_1^{-1}(W_1(\Omega))=\binom{n}{2k}(2k)! \omega_n \(\frac{|M|}{\omega_n}\)^{\frac{n+1-2k}{n}}.
	\end{align}
	Equality holds in \eqref{s9:Lk-shifted-quermassintegral} or \eqref{s9:Lk-area} if and only if $M$ is a geodesic sphere centered at the origin.	
\end{proof}

\end{document}